\documentclass[a4paper,1p,sort&compress,dvips,10pt]{elsarticle}
\usepackage{latexsym}
\usepackage{amsmath}
\usepackage{amsthm}
\usepackage{amsfonts}
\usepackage{amssymb}
\usepackage{bm}
\usepackage[english]{babel}

\DeclareFontFamily{OT1}{pzc}{}
\DeclareFontShape{OT1}{pzc}{m}{it}%
             {<-> s * [1.125] pzcmi7t}{}
\DeclareMathAlphabet{\mathpzc}{OT1}{pzc}{m}{it}

\journal{{\tt arXiv}, typeset with elsarticle.cls}

\newcommand{\be}{\begin{equation}}
\newcommand{\ee}{\end{equation}}

\newcommand{\ef}[1]{\, #1}

\newcommand{\cO}{{\cal O}}

\renewcommand{\emptyset}{\varnothing}
\newcommand{\smfrac}[2]{\genfrac{}{}{0.25pt}{1}{#1}{#2}}

\newcommand{\modd}[1]{\ \mathrm{(mod~}#1\mathrm{)} }

\newcommand{\fpl}{\mathpzc{Fpl}\!}
\newcommand{\LP}{\mathpzc{LP}}

\newcommand{\Sym}{\mathrm{Sym}}

\newcommand{\proofof}[1]{\noindent {\it Proof of #1:~}}

\newcommand{\sqqV}{\,\rule{.5pt}{7pt}\rule{6pt}{0pt}\rule{.5pt}{7pt}\,}
\newcommand{\sqqH}{\,\makebox[0pt][l]{\rule{7pt}{.5pt}}\raisebox{6.5pt}{\rule{7pt}{.5pt}}\,}

\newcommand{\beo}{{\bf e}}
\newcommand{\buno}{{\bf 1}}

\newcommand{\simX}{\includegraphics{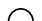}}
\newcommand{\nsimX}{\includegraphics{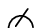}}

\newcommand{\kket}[1]{\mbox{$\| #1 \rangle \! \rangle $}}

\newcommand{\ket}[1]{\mbox{$| #1 \rangle $}}
\newcommand{\bra}[1]{\mbox{$\langle #1 | $}}
\newcommand{\braket}[2]{\big\langle #1 | #2 \big\rangle}

\newtheorem{theorem}{Theorem}[section]
\newtheorem{corollary}{Corollary}[section]

\newtheorem{lemma}{Lemma}[section]

\newtheorem{proposition}{Proposition}[section]

\begin{document}

\begin{frontmatter}

\title
{A one-parameter refinement\\of the Razumov--Stroganov correspondence}

\author[cer]{Luigi Cantini}
\ead{lu{}igi.can{}tini@u-ce{}rgy.fr}
\author[mi]{Andrea Sportiello}
\ead{Andr{}ea.Sportie{}llo@mi.in{}fn.it}

\address[cer]{Universit\'e de Cergy-Pontoise, LPTM -- UMR 8089 of
  CNRS\\
2 av.~Adolphe Chauvin, 95302 Cergy-Pontoise, France}
\address[mi]{Dipartimento di Fisica dell'Universit\`a degli Studi di
  Milano, and INFN,\\
  via Giovanni Celoria 16, 20133 Milano, Italy}

\date{February 23xx, 2012}

\begin{abstract}
We introduce and prove a one-parameter refinement of the
Razumov--Stroganov correspondence. This is achieved for fully-packed
loop configurations (FPL) on domains which generalize the square
domain, and which are endowed with the gyration operation. We consider
one given side of the domain, and FPLs such that the only
straight-line tile on this side is black. We show that the enumeration
vector associated to such FPLs, weighted according to the position of
the straight line and refined according to the link pattern for the
black boundary points, is the ground state of the \emph{scattering
  matrix}, an integrable one-parameter deformation of the $O(1)$
Dense Loop Model Hamiltonian.  We show how the original
Razumov--Stroganov correspondence, and a conjecture formulated by Di
Francesco in 2004, follow from our results.
\end{abstract}


\end{frontmatter}

\section{Introduction}

\noindent
The \emph{Razumov--Stroganov correspondence}
\cite{rs,usRS} relates some fine statistical properties of two
distinct integrable systems in Statistical Mechanics \cite{baxter}: on
one side the
\emph{6-Vertex Model} on portions of the square lattice, with
domain-wall boundary conditions, and on the other side the
\emph{$O(1)$ Dense Loop Model} (DLM) with cyclic boundary
conditions.

The configurations of the first model have several easy
reformulations, in terms of \emph{Fully-packed loops} (FPL),
\emph{Alternating Sign Matrices}, or a family of 
monotone arrays
called \emph{Gog triangles} \cite{bress}. In the FPL
incarnation, to each configuration $\phi$ is naturally associated a
non-crossing pairing $\pi = \pi(\phi)$ of a set of cyclically-ordered
points (called \emph{link pattern}). We call $\Psi_{\rm FPL}(\pi)$ the
corresponding enumerations, i.e.\ the number of $\phi$'s such that
$\pi(\phi)=\pi$, and $Z_{\rm FPL} = \sum_{\pi} \Psi_{\rm FPL}(\pi)$
their total number. An explicit formula is known for this quantity
\cite{zeil,kupe},
\be
Z_{\rm FPL}(n) 
= \prod_{j=0}^{n-1} \frac{(3j+1)!}{(n+j)!}  
\ef.  
\ee 
For the second model, we have reformulations in terms of the
\emph{integrable XXZ quantum spin chain} at $\Delta = -\frac{1}{2}$,
and the \emph{Potts Model} at the percolation point, $Q=1$. The model
is realised on a semi-infinite cylinder, and is naturally analysed
through transfer-matrix techniques. In the DLM incarnation, the
transfer matrix $T_{\pi\, \pi'}$ acts on a space whose states are
naturally labeled by link patterns. The matrix $T$ encodes the
transition rates of a Markov Chain on this space, and there is a
unique steady-state distribution $\Psi_{O(1)}(\pi)$, called
\emph{ground state}, and corresponding to the Frobenius right
eigenvector of $T$. Integrability shows that $T$ commutes with a
simpler operator, the \emph{Hamiltonian}, $H_0$, so that
$\Psi_{O(1)}(\pi)$ is also a right eigenvector of $H_0$.  As the
corresponding left eigenvector is the uniform vector, with all entries
equal to 1, the natural norm of $\Psi_{O(1)}(\pi)$ is given by the sum
of the entries, $Z_{O(1)} = \sum_{\pi} \Psi_{O(1)}(\pi)$.

The Razumov--Stroganov correspondence states that, under the
normalisation for $\Psi_{O(1)}$ that sets $Z_{O(1)} = Z_{\rm FPL}$, we
have $\Psi_{\rm FPL}(\pi) = \Psi_{O(1)}(\pi)$ for all link
patterns~$\pi$. This fact was conjectured in \cite{rs}, and proven by
the authors in~\cite{usRS}.

A great effort has been devoted to the study of the properties of the
ground state of the $O(1)$ Dense Loop Model.  Building on the
integrable structure of the DLM, some deep connections with the
representation theory of $U_q(\widehat{sl_2})$, or of Affine Hecke
Algebras, have been elucidated, and even connections with algebraic
geometry have emerged (see \cite{PZJthesis} for a review).

The power of integrability manifests itself when the original loop
model is deformed introducing the so-called \emph{spectral parameters}
$\vec{z} = \{z_i\}$ (the uniform counting corresponds to the choice
$z_i=1$ for all $i$). The components of the ground state,
$\Psi_{O(1)}(\pi)$, that, in the uniform model and after
normalisation, are all integers, are deformed into polynomials in
these parameters.  Besides the emergence of the connections mentioned
above, this procedure has more concretely allowed to obtain closed
formulas for certain linear combinations of components of the ground
state, through determinantal
representations, or through multiple contour
integral formulas (see \cite{PZJthesis} and references therein).  In
particular, the deformation of the normalisation of the $O(1)$ ground
state, $Z_{O(1)}(\vec{z})$, can be evaluated, if the value at a
reference pattern $\pi$ of `rainbow' shape is fixed~\cite{pdf-pzj1}.

Analogously, on the FPL side (in its equivalent formulation as
6-Vertex Model at \hbox{$\Delta=1/2$}), following the general strategy
of Yang--Baxter integrability, it is quite natural to introduce
spectral parameters, associated to row- and column-indices of the
square grid, and deform both the total enumeration of configurations,
$Z_{\rm FPL}$, and the refined enumerations, $\Psi_{\rm FPL}(\pi)$,
into polynomials in these variables.  The polynomial 
$Z_{\rm FPL}(\vec{z})$ has a determinantal representation, called
\emph{Izergin-Korepin determinant} \cite{IK_I}. When specialised at
$q=e^{\frac{2 \pi i}{3}}$, again with a natural normalisation of the
`rainbow' reference patterns, it coincides with its DLM analogue
$Z_{O(1)}(\vec{z})$~\cite{pdf-pzj1}.

A natural question one could pose is whether and to which extent it is
possible to introduce parameters both on the DLM and FPL sides, in
such a way to produce deformations of the Razumov--Stroganov
correspondence, i.e.\ polynomial identities at the level of the
refined enumerations $\Psi_{\rm FPL}(\pi)$ and $\Psi_{O(1)}(\pi)$.

An attempt in this direction was pursued by Di Francesco
\cite{PdF04}. The proposal consisted, on the FPL side, in taking all
spectral parameters equal to 1, except for the one associated to the
bottom row of the square, valued $z$.  We say that a FPL configuration
has \emph{refinement position} $j$ if, on the bottom row, the unique
tile consisting of a straight line is at column $j$ (in the ASM
representation, this corresponds to say that this is the position of
the unique $+1$ in the bottom row).  The total number $Z_{\rm
  FPL}^{[j]}(n)$ of FPL's with refinement position $j$ is also known,
and given by the formula \cite{zeilRef}
\be
Z_{\rm FPL}^{[j-1]}(n-1)
=
\frac{(n+j)!\, (2n-j)!\, (2n+1)!}{n!\, j!\, (n-j)!\, (3n+1)!}
\; Z_{\rm FPL}(n-1)
\ef.
\ee
%
The introduction of the spectral parameter $z$ corresponds to count
with a weight $t^{j-1}$ each configuration having refinement position
$j$, for $t=t(z)=\frac{qz-1/q}{q-z/q}$.  Thus, on the FPL side, we
have counting polynomials $\Psi_{\rm FPL}(t;\pi)$, that reduce to
$\Psi_{\rm FPL}(\pi)$ for $t=1$.

Based on numerical experimentations, Di Francesco conjectured that,
while the set of $\Psi_{\rm FPL}(t;\pi)$'s does not match the ground
state of any known integrable deformation of the $O(1)$ loop model,
its symmetrisation under rotation is equal to the symmetrisation of
$\Psi^{(i)}_{O(1)}(t;\pi)$, the unique ground state of the
\emph{scattering matrix at site $i$}, $S_i(t)$ (for any $1 \leq i \leq
2n$), which is a one--parameter deformation of the Hamiltonian of the
$O(1)$ loop model. At the light of a `dihedral covariance' of the
ground states $\Psi^{(i)}_{O(1)}(t;\pi)$ (discussed in detail later
on), one can concentrate on the state $\Psi^{(1)}_{O(1)}(t;\pi)$.


In the present paper we address and prove Di Francesco's conjecture by
actually proving a stronger refined correspondence, that does not
require a symmetrisation.

As a first direction of generalisation, we shall consider FPL's not
only on a regular square, but on a family of domains that we shall
call \emph{dihedral domains}. These domains are characterized by the
existence of two \emph{gyration operations}, implying the invariance
under rotation of the usual FPL enumerations, in a way that extends
the original work of Wieland \cite{wie}. This family of domains
includes some of the ``symmetry classes'' of FPL (or equivalently of
Alternating Sign Matrices) for which a Razumov--Stroganov conjecture
was formulated, namely HTASM and QTASM (half-turn and quarter-turn
symmetric) \cite{Raz-Str-2}, but it is actually much larger.  This
extension was in fact already presented in \cite{usRS}, and the expert
reader should not be surprised by the fact that this family of domains
is still the appropriate setting also in the framework of the Di
Francesco's conjecture mentioned above.

Most importantly, we introduce and study FPL's $\phi$, enumerated
according to link patterns which are associated to $\phi$'s through a
function $\tilde{\pi}(\phi)$ which is \emph{different} from the one
usually considered in the literature. This gives new enumerations
$\tilde{\Psi}_{\rm FPL}(\pi)$ and corresponding refined enumerations
$\tilde{\Psi}_{\rm FPL}(t;\pi)$. On one side, we prove that
$\tilde{\Psi}_{\rm FPL}(t;\pi)$ and $\Psi^{(1)}_{O(1)}(t;\pi)$ coincide,
with no need of symmetrisation. This is the result of the paper that
we consider structurally more relevant.
On the other side, we prove that the symmetrisation of
$\tilde{\Psi}_{\rm FPL}(t;\pi)$ coincides with the symmetrisation of
$\Psi_{\rm FPL}(t;\pi)$.

For concreteness, in this introduction we describe the function
$\tilde{\pi}(\phi)$ when the domain $\Lambda$ is the $n \times n$
square (the case of general domains is treated in depth in the body of
the paper).  We recall that a FPL is a bicolouration (in black and
white) of the edges of the domain $\Lambda$, such that each internal
vertex is adjacent to two black and two white edges.  Let
$\fpl(\Lambda)$ be the ensemble of FPL on $\Lambda$ in which the
external edges are coloured alternatively black and white.  Thus, the
colouration of a single reference external edge, say the vertical one
at the bottom left corner, completely determines the boundary
conditions, and $\fpl(\Lambda)$ is the disjoint union of the two sets
$\fpl_+(\Lambda)$ and $\fpl_-(\Lambda)$. The map $\sigma$, consisting
in swapping black and white, is an involution on $\fpl(\Lambda)$, and
a bijection between $\fpl_+(\Lambda)$ and $\fpl_-(\Lambda)$.

In the literature, when referring to the ``FPL side'' of the
Razumov--Stroganov correspondence, it was always meant that FPL's were
in the ensemble $\fpl_+(\Lambda)$ (or $\fpl_-(\Lambda)$), and were
refined according to the black link pattern, with indices from 1 to
$2n$ assigned to the black external edges once and for all.  The
function $\pi(\phi)$, that associates a link pattern $\pi$ to a
configuration $\phi \in \fpl_+(\Lambda)$, crucial in the definition of
$\Psi_{\rm FPL}(\pi)$ and thus of the Razumov--Stroganov
correspondence, is the one given by such a prescription.

Wieland gyration implies as a corollary that the refined enumerations
on $\fpl_+(\Lambda)$ and on $\fpl_-(\Lambda)$ coincide. Of course, the
involution $\sigma$ exchanges the black and white link patterns
associated to a configuration, so the statement of the
Razumov--Stroganov correspondence holds as well for the white link
patterns. 

Consider a generic class of functions $\pi'(\phi)$, that associate to
$\phi$ the black \emph{or} the white link pattern, depending on some
properties of $\phi$, and with the black (or white) external edges
numbered from $1$ to $2n$, in counter-clockwise order, starting from
some external edge \emph{depending from $\phi$}. In other words,
$\pi'(\phi)$ is determined completely by a choice of external edge
$e(\phi)$ for $\phi$ (by setting the colour of the link pattern to the
colour of edge $e$, and assigning the label $1$ to $e$. As a
consequence, $\pi'(\phi) \equiv \pi'(\sigma \phi)$).  The function
$\pi(\phi)$ described above is such that the function $e(\phi)$ is
\emph{constant}, i.e.\ a reference external edge is fixed once and for
all. The enumeration is restricted to $\fpl_+(\Lambda)$.

Even at the light of the ordinary Razumov--Stroganov correspondence,
except for the choice of a constant function for $e(\phi)$, there is
no special reason \emph{a priori} for hoping that the enumerations
induced by functions $\pi'(\phi)$ have any remarkable property, and in
particular any relation with $\Psi^{(i)}_{O(1)}(t;\pi)$.

\begin{figure}
\begin{center}
\makebox[0pt][c]{$\pi(\phi)$}
\rule{281.6pt}{0pt}
{}
\makebox[0pt][c]{$\tilde{\pi}(\phi)$}

\vspace{-1mm}
\includegraphics[scale=0.8]{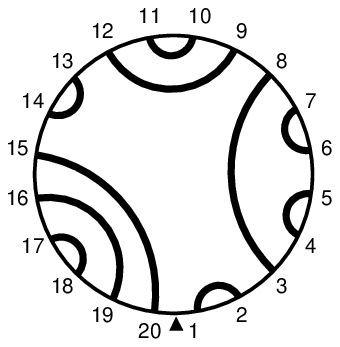}
\includegraphics[scale=1.6]{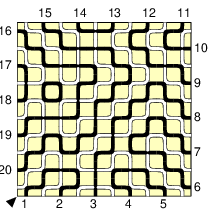}
\includegraphics[scale=1.6]{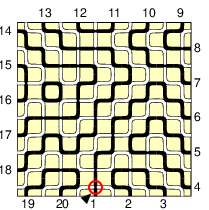}
\includegraphics[scale=0.8]{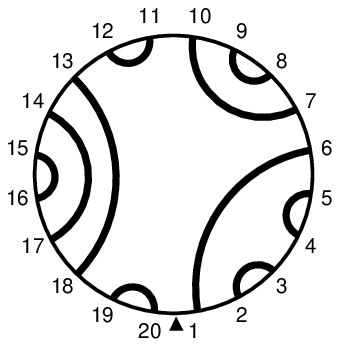}

\vspace{2mm}
\includegraphics[scale=0.8]{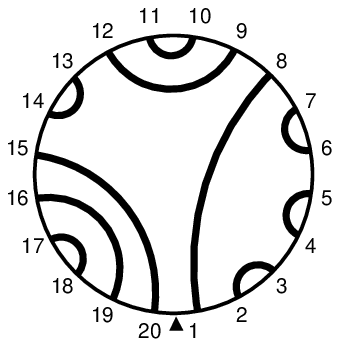}
\includegraphics[scale=1.6]{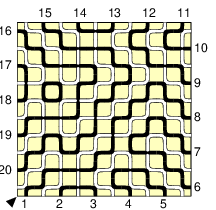}
\includegraphics[scale=1.6]{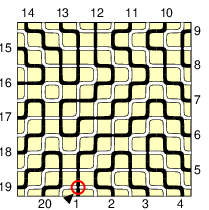}
\includegraphics[scale=0.8]{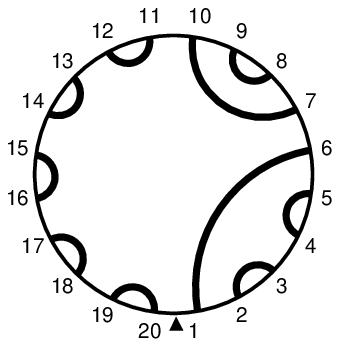}
\end{center}
\caption{\label{fig.expitilde}Comparison between the `ordinary' map
  $\pi(\phi)$, on the left, and the `new' map $\tilde{\pi}(\phi)$, on
  the right. Two distinct examples are shown in the two lines. In the
  first one, the refinement position is odd, $h(\phi)=5$, thus
  $\tilde{\pi}(\phi)$ is the appropriate rotation of the black link
  pattern of $\phi$. In the second line, $h(\phi)=4$, even, thus
  $\tilde{\pi}(\phi)$ is the appropriate rotation of the \emph{white}
  link pattern of $\phi$.}
\end{figure}

Nonetheless, the refinement on the bottom row suggests a different
natural, non-uniform choice for $e(\phi)$: we take as $e(\phi)$ the
external edge incident to the refinement position. We call
$\tilde{\pi}(\phi)$ the function $\pi'(\phi)$ corresponding to this
choice.  The enumerations $\tilde{\Psi}_{\rm FPL}(t;\pi)$, with the
properties anticipated above, are the ones obtained by using this
function $\tilde{\pi}(\phi)$, and associating a weight $t^{j-1}$ to a
configuration with refinement position $j$.  The enumeration is
restricted to the set $\fpl_b(\Lambda)$, of FPL's such that the
external edge incident to the refinement position is black.
One of our main results, Theorem \ref{main}, states that, for a
general class of domains (including the $n \times n$ square as a
special case), the vector $\tilde{\Psi}_{\rm FPL}(t;\pi)$ is an
eigenvector of the scattering matrix $S_1(t)$, and thus is
proportional to $\Psi_{O(1)}^{(1)}(t;\pi)$.
A comparison between the ordinary $\pi(\phi)$ and this new map
$\tilde{\pi}(\phi)$ is also shown through examples in
Figure~\ref{fig.expitilde}.

We will then show how Di Francesco's conjecture, now reducing to
$\Sym\,\tilde{\Psi}_{\rm FPL}(t;\pi) = \Sym\,\Psi_{\rm FPL}(t;\pi)$, 
where $\Sym$ is a symmetrisation operator, follows from a study of the
orbits under the action of Wieland half-gyration, and the associated
evolution of the refinement position. Again, this is done for the
general class of domains.

As already shown in \cite{PdF04}, Di Francesco's conjecture is a
generalisation of the ordinary Razumov--Stroganov correspondence, thus
this analysis also provides an alternative proof for the
latter. However, there are various other alternative proofs that can
be derived at various stages of the analysis.  In particular, one
derivation involves our sole Theorem \ref{main}, plus a non trivial
bijection between $\fpl_{b}(\Lambda)$ and $\fpl_{+}(\Lambda)$, which
preserves the link pattern associated to a configuration.


\bigskip
\noindent
The paper is organized as follows.  In Section \ref{Temperley-section}
we introduce the Cyclic Temperley--Lieb Algebra, which is at the basis
of the integrable structure of the $O(1)$ Dense Loop Model. In Section
\ref{exch-section} we define the Temperley--Lieb Hamiltonian and the
scattering equations, and derive some relevant properties of their
solutions.  In Section \ref{sec.FPLdefs} we begin by recalling some
basic facts about FPL's, then we proceed with the definition of
dihedral domains. In Section \ref{sec.WieRemind}, we define the
gyration operations acting on FPL, and show that they act as a
rotation at the level of the link pattern associated to the FPL. This
is the generalisation of the Wieland gyration theorem to our family of
domains, and is essentially a reminder of the theory already presented
in \cite{usRS}.
Sections \ref{main-section} and \ref{from-main-to-DF} produce the results
sketched above, in particular Theorem \ref{main}, stating that the
`new' FPL enumeration provides a solution of the scattering equation,
and Theorem \ref{df-conj}, the Di Francesco's 2004 former conjecture,
relating the `old' and `new' FPL enumerations.



\section{The Cyclic Temperley--Lieb Algebra and the $O(1)$ Dense Loop
  Model}
\label{Temperley-section} 

\noindent
In this section we analyse the $O(1)$ Dense Loop Model side of the
correpondence, which consists of the Perron--Frobenius eigenvector
associated to the \emph{scattering equation}, a linear equation
involving a representation of the Temperley--Lieb Algebra.  In Section
\ref{ssec.TLdef} this algebra is defined, while in Section
\ref{exch-section} we define our vector of interest, and deduce some
of its properties.

\subsection{The Cyclic Temperley--Lieb Algebra}
\label{ssec.TLdef}

\noindent
We start by recalling the definiton of the 
\emph{Cyclic Temperley--Lieb Algebra} $\mathrm{CTL}_{N}(\tau)$, which
is the free algebra with generators $\{\beo_i \}_{i \in \mathbb{Z}}$,
and the invertible rotation operator $R$, and relations
\begin{subequations}
\begin{align}
\beo_i
&=
\beo_{i+N}
\ef;
&
R^{\pm 1} \beo_i 
&=
\beo_{i \pm 1} R^{\pm 1}
\ef;
\label{eq.6537665}
\\ 
\beo_i^2
&=
\tau \beo_i
\ef;
&
\beo_i \beo_{i\pm 1} \beo_i
&=\beo_i
\ef;
\\
[\beo_i,\beo_j]
&=0
\makebox[0pt][l]{\rule{30pt}{0pt}
\textrm{for $i-j \not\equiv \pm 1 \ \modd{N}$.}}
&& {}
\end{align}
\end{subequations}
This algebra has interesting properties for a full range of the
parameter $\tau$, and remarkable specialisation at a family of
discrete values for $\tau$ (an useful alternative parametrisation is
to set $\tau=-q-q^{-1}$). In the following we shall restrict to
$\tau=1$, i.e.\ $q$ a cubic root of unity, and consider two
kinds of diagrammatic representations of $\mathrm{CTL}_{N}(\tau)$,
acting on spaces of \emph{link patterns}.

\begin{figure}[!t]
\begin{center}
\includegraphics[scale=0.55]{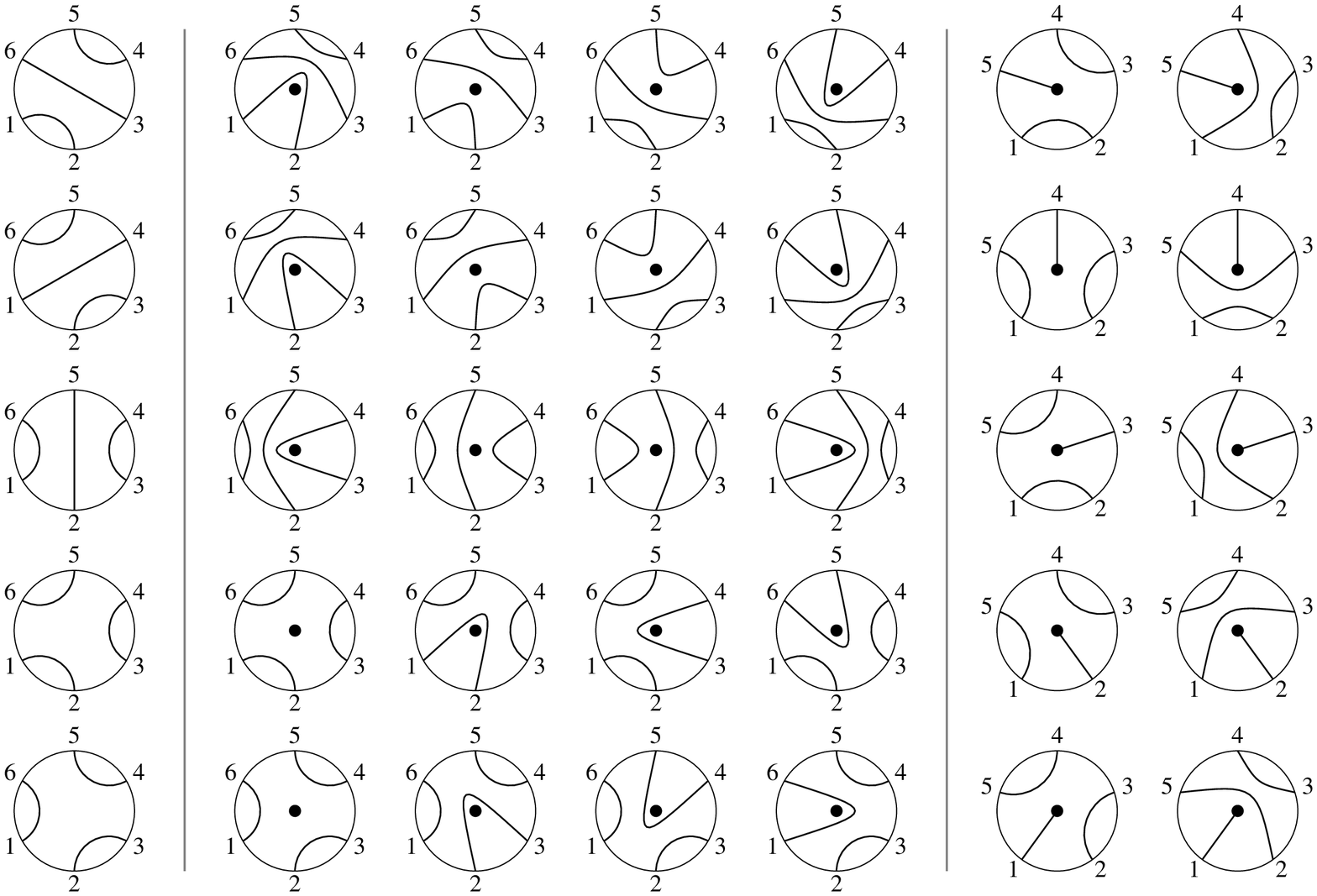}
\end{center}
\caption{\label{fig.LPvari}The sets $\LP(6)$, $\LP^*(6)$ and
  $\LP^*(5)$.}
\end{figure}

Define $\LP(2n)$ as the set whose elements are \emph{link patterns}
with $n$ arcs, i.e., the possible topologies of $n$ arcs that connect,
through a non-crossing pairing, $2n$ points ordered cyclically
counter-clockwise along the boundary of a disk.  The set $\LP^*(N)$ is
an analogue of $\LP(2n)$ consisting of \emph{punctured link
  patterns}. For $N=2n$ even, these are the possible topologies of $n$
arcs that connect, through a non-crossing pairing, $2n$ points ordered
cyclically along the boundary of a punctured disk (thus, it matters if
an arc passes on the right or on the left of the punture).  For
$N=2n-1$ odd, these are the possible topologies of $n$ arcs that
connect, through a non-crossing pairing, $2n-1$ points ordered
cyclically along the boundary of a punctured disk, and the punture
itself.

It is easily seen that the sets $\LP(2n)$, $\LP^*(2n)$ and
$\LP^*(2n+1)$ have cardinalities $C_n$, $(n+1) C_n$ and $(2n+1) C_n$,
where $C_n = \frac{1}{n+1} \binom{2n}{n}$ is the $n$-th Catalan
number.  See Figure \ref{fig.LPvari} for an illustration of these
patterns.

It is useful to establish a reference pattern. We define a
\emph{rainbow pattern} as a link pattern consisting of a unique
`rainbow' of parallel arcs
$\big\{\ldots,\{i,j\},\{i+1,j-1\},\{i+2,j-2\},\ldots\big\}$, with the
possible puncture maximally nested. There are $n$ rainbow patterns in
$\LP(2n)$, and $N$ ones in $\LP^*(N)$, related by rotation.  Examples
for $\LP(8)$, $\LP^*(8)$ and $\LP^*(7)$ are 
\be
\label{eq.Rainbdef}
\raisebox{-20pt}{\includegraphics[scale=0.55]{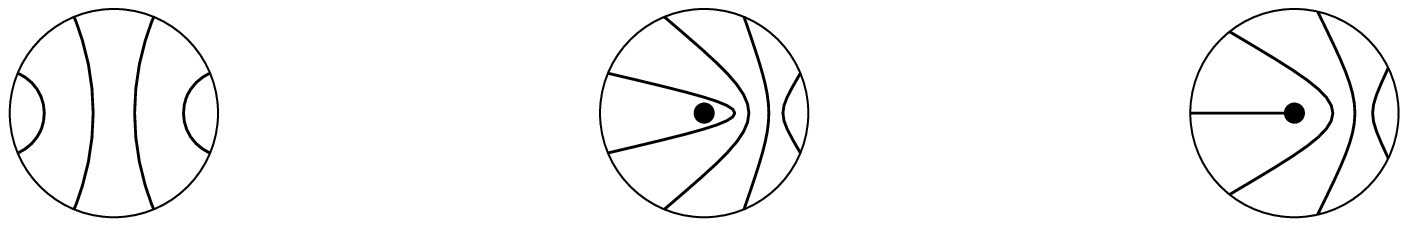}}
\ee
The first representation of $\mathrm{CTL}_{N}(\tau)$ we will consider
is defined for $N=2n$ even, and acts on the vector space
$\mathbb{C}^{\LP(2n)}$, having a priviliged basis whose elements are
indicised by elements of $\LP(2n)$. The second class of
representations acts on $\mathbb{C}^{\LP^*(N)}$, with basis elements
indicised by elements of $\LP^*(N)$.

The action on the basis vectors is induced by the graphical action of
the Temperley--Lieb generators on the link patterns.  The operators
$\beo_j$ and $R$\;\footnote{Note that equation (\ref{eq.6537665}), and
  the consistent representation in (\ref{eq.defEandR}), fix a direction
  convention for $R$, which rotates the link pattern by
  \emph{increasing} the indices, $R\big( \big\{ \ldots,\{i,j\},\ldots
  \big\} \big) = \big( \big\{ \ldots,\{i+1,j+1\},\ldots \big\}
  \big)$.}  are represented as customarily as
\begin{align}
\label{eq.defEandR}
\beo_j:
\quad
&
\setlength{\unitlength}{11pt}
\raisebox{-32pt}{
\begin{picture}(6.5,5.5)
\put(0,0){\includegraphics[scale=0.55]{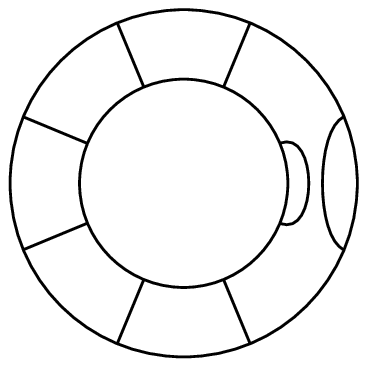}}
\put(0,4.2){$\scriptstyle{N}$}
\put(0.2,2){$\scriptstyle{1}$}
\put(1.8,0.3){$\scriptstyle{2}$}
\put(2.6,0.2){$\scriptstyle{\ldots}$}
\put(2,3.6){$\scriptstyle{N}$}
\put(2.1,2.6){$\scriptstyle{1}$}
\put(2.65,2.15){$\scriptstyle{2\ldots}$}
\put(2.6,0.2){$\scriptstyle{\ldots}$}
\put(3.3,3.6){$\scriptstyle{j+1}$}
\put(4.1,2.6){$\scriptstyle{j}$}
\put(5.9,4.2){$\scriptstyle{j+1}$}
\put(5.9,2){$\scriptstyle{j}$}
\end{picture}}
&
R:
\quad
&
\setlength{\unitlength}{11pt}
\raisebox{-32pt}{
\begin{picture}(6.5,2)
\put(0,0){\includegraphics[scale=0.55]{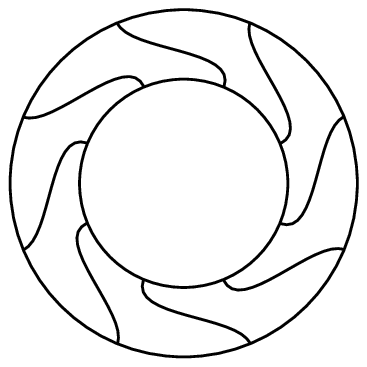}}
\put(0,4.2){$\scriptstyle{N}$}
\put(0.2,2){$\scriptstyle{1}$}
\put(1.8,0.3){$\scriptstyle{2}$}
\put(2.6,0.2){$\scriptstyle{\ldots}$}
\put(2,3.6){$\scriptstyle{N}$}
\put(2.1,2.6){$\scriptstyle{1}$}
\put(2.65,2.15){$\scriptstyle{2\ldots}$}
\put(2.6,0.2){$\scriptstyle{\ldots}$}
\end{picture}}
\end{align}
Then, the action corresponds to the concatenation of the diagrams (an
annulus, corresponding to an operator, is juxtaposed outside a disk,
containing a link pattern, and a new pattern is derived from the
topology of the matching of the new boundary points). In case loops
are produced, they are dropped, with a factor $\tau$ each (i.e., in
our specialisation to $\tau=1$, we just forget about loops).

It is easily seen that the representation on $\LP(2n)$ is a quotient
of the representation on $\LP^*(2n)$, obtained by just `forgetting'
about the location of the puncture.

For two terminations $i, j \in \{1,\ldots,N\}$, and a link pattern
$\pi$, we write that $i \sim j$ if $(ij)$ is an arc of $\pi$, and $i
\nsim j$ otherwise. If the terminations are consecutive, we say that
$i \simX i+1$ if $\pi$ is in the image of $e_i$, and $i \nsimX i+1$
otherwise. In $\LP(2n)$ and $\LP^*(2n+1)$, $i \simX i+1$ if and only
if $i \sim i+1$, while in $\LP^*(2n)$, $i \simX i+1$ if and only if $i
\sim i+1$ and this short arc does not embrace the puncture.

Working in the linear spaces $\mathbb{C}^{\LP(2n)}$ and
$\mathbb{C}^{\LP^*(N)}$ we use a ``ket'' notation for vectors
$\ket{v}=\sum_\pi v(\pi) \ket{\pi} $, with the sum over $\pi$ running
in the appropriate ensemble $\LP(2n)$ or $\LP^*(N)$ depending on
circumstances. We denote covectors as $\bra{v}$, and scalar products
$\sum_\pi u(\pi)^* v(\pi)$ as $\braket{u}{v}$ (in fact, complex
conjugation is never used here, and we could have worked equivalently
in the linear spaces $\mathbb{R}^{\LP(2n)}$ or
$\mathbb{Q}^{\LP(2n)}$). Such a notation is chosen in order to
distinguish easily vectors, linear operators and covectors, while
omitting the specification of the space, between
$\mathbb{C}^{\LP(2n)}$ and $\mathbb{C}^{\LP^*(N)}$, writing single
equations that hold simultaneously in the two cases.  The maps
$\beo_i$ and $R$ induce linear operators on $\mathbb{C}^{\LP(2n)}$ and
$\mathbb{C}^{\LP^*(N)}$, for which, with abuse of notation, we
continue to use the same symbols.

\subsection{The Hamiltonian, the scattering matrix and the scattering
  equations} 
\label{exch-section}

\noindent
Consider a special linear expression in Temperley--Lieb Algebra,
named \emph{Hamiltonian}
\be\label{hamil-O(1)}
H_0=\sum_{i=1}^{N}(\beo_i-\buno)
\ef.
\ee
The action of this operator, both on $\mathbb{C}^{\LP(2n)}$ (if $N=2n$) and on
$\mathbb{C}^{\LP^*(N)}$, has at sight a unique left null-vector,
$\bra{1}=\sum_\pi \bra{\pi}$. The \emph{ground state} of the Hamiltonian
is the unique right null-vector, i.e.\ the unique vector $\ket{\Psi_{O(1)}}
\in \mathbb{C}^{\LP(2n)}$ (or in $\mathbb{C}^{\LP^*(N)}$) satisfying
\be
\label{RS-eq0}
H_0 \ket{\Psi_{O(1)}} =0
\ef.
\ee
In this section we omit the $O(1)$ subscript for brevity.  The vector
$\ket{\Psi}$ can clearly be normalised in such a way that all of its
entries are integers (because $H_0$ has integer entries). Furthermore,
all the entries are positive. This comes from the stochastic nature of
$H_0$ (i.e., $\ket{\Psi}$ is the Perron--Frobenius vector for the
matrix $\frac{1}{N} H_0+ \buno$).

Following closely \cite{PdF04}, we now introduce a family of
equations, called \emph{scattering equations}, depending on one
parameter $t$, that generalise (\ref{RS-eq0}), in the sense that the
solutions $\ket{\Psi^{(i)}(t)}$ of the $i$-th equation at parameter
$t$ reduce to the ground state $\ket{\Psi}$ of the Hamiltonian, when
$t=1$.  Introduce the operator $X_i(t)$\;\footnote{For the educated
  reader, this operator is related to the usual `baxterization' of the
  Cyclic Temperley--Lieb algebra, through
\[
\check{\mathrm{R}}_i(z)= X_i\left(\smfrac{qz-q^{-1}}{q-q^{-1}z}\right)
\ef,
\]
  setting $q=e^{\frac{2\pi i}{3}}$.}
\be
X_i(t) = t\; \buno + (1-t)\, \beo_i  
\ef.
\ee
The \emph{scattering equation} we are going to consider is 
\be
\label{eq.qKZfund}
X_i(t)\ket{\Psi^{(i)}(t)}= R\,\ket{\Psi^{(i)}(t)}
\ef.
\ee
%
A simple analysis of equation (\ref{eq.qKZfund}) tells us that its
solution is unique up to normalization.
Indeed, if $0<t<1$, the matrix $R\,X_i(t)$ is a Markov chain and it is
not difficult to show that it is irreducible.  Therefore existence and
unicity (up to normalization) of the solution of equation
(\ref{eq.qKZfund}) follow from the Perron--Frobenius theorem, in a way
completely analogous to the reasoning done for $\ket{\Psi}$ and the
Hamiltonian. The components $\Psi^{(i)}(t;\pi)$ of the vector
$\ket{\Psi^{(i)}(t)}$ can be normalized in such a way that they are
polynomials of the parameter $t$ and the Perron--Frobenius theorem
ensures their positivity, when evaluated on the open real interval $t
\in \;]0,1[$, but remarkably \emph{not} the observed, and proven later
on, coefficient-wise positivity.

In general $\ket{\Psi^{(i)}(t)}$ is not dihedrally invariant (as we
will see later on, it is only ``dihedrally covariant'', for a
simultaneous action on link patterns $\pi$ and superscript $(i)$).
The dihedral symmetry is however restored at $t=1$, where
$\ket{\Psi^{(i)}(t)}$ reduces to the solution of equation
(\ref{RS-eq0}). Indeed, for $t \to 1$, the operators $X_i(t)$ reduce
to the identity, and therefore from equation (\ref{eq.qKZfund}) it
follows immediately that $\ket{\Psi^{(i)}(1)}$ is rotationally
invariant
\be
\label{rot-inv}
\ket{\Psi^{(i)}(1)}=R\,\ket{\Psi^{(i)}(1)}
\ef.
\ee
Furthermore, if we apply the projector
\be
\Sym := \frac{1}{N}\sum_{i=0}^{N-1} R^k
\ee
on both sides of equation (\ref{eq.qKZfund}) (and use the obvious
$\Sym\, R=\Sym$) we find
\be
\label{sym(e-1)}
\Sym\,(\beo_i-\buno)\ket{\Psi^{(i)}(t)}=0
\ef.
\ee
Combining equations (\ref{rot-inv}) and (\ref{sym(e-1)}), we obtain
\be
\begin{split}
H_0\ket{\Psi^{(i)}(1)} 
&= 
\sum_{k=0}^{N-1}
R^k(\beo_{i}-\buno)R^{-k} \ket{\Psi^{(i)}(1)} 
= \sum_{k=0}^{N-1}
R^k(\beo_{i}-\buno) \ket{\Psi^{(i)}(1)} 
\\
&= N\, \Sym\,
(\beo_i-\buno)\ket{\Psi^{(i)}(1)} =0 
\ef;
\end{split}
\ee
i.e., $\ket{\Psi^{(i)}(1)}$ is proportional to the ground state of the
Hamiltonian.  An alternate derivation of the same result, already
presented in \cite{PdF04}, goes through the definition of the
\emph{scattering matrix}
\be
S_i(t)=X_{i+N-1}(t) \cdots X_{i+1}(t)X_i(t)
\ef.
\ee
Since $R^{k} X_i(t)R^{-k} = X_{i+k}(t)$, we find that
\be
X_{i+k}(t) X_{i+k-1}(t) \dots X_{i+1}(t) X_i(t)
=
R^{k+1}(R^{-1}X_i(t))^{k+1}
\ef,
\ee
and therefore, if $\Psi^{(i)}(t)$ satisfies equation (\ref{eq.qKZfund}),
then it must satisfy 
\be
X_{i+k}(t) X_{i+k-1}(t) \dots X_{i+1}(t) X_i(t)
\ket{\Psi^{(i)}(t)}
= R^{-k-1} \ket{\Psi^{(i)}(t)}
\ef,
\ee
and in particular
\be
\label{scatter-eq}
\big( S_i(t) - \buno \big)
\ket{\Psi^{(i)}(t)}= 0
\ef.
\ee
Observe that, for all $i$, the expansion of $S_i(t)$ near $t=1$ leads to
the Hamiltonian. Indeed
$S_i(1) = \buno$, and
\begin{align}
\left.
\frac{dS_i(t)}{dt}
\right|_{t=1} 
&=
\sum_{j=i}^{i+N-1}
\left.
X_{i+N-1}(t)
\cdots
X_{j+1}(t)
\frac{dX_{j}(t)}{dt}
X_{j-1}(t)
\cdots
X_{i}(t)
\right|_{t=1}
=
-H_0
\ef;
\end{align}
therefore, taking the derivative of equation (\ref{scatter-eq}) and then
setting $t=1$, one finds that $\Psi^{(i)}(1)$ satisfies equation
(\ref{RS-eq0}).

As anticipated, the solutions of 
equation (\ref{eq.qKZfund}) show
\emph{dihedral covariance} 
for different values of
$i$. 
Indeed we find that 
\be
0 = \big( X_i(t)-R \big)\,\ket{\Psi^{(i)}(t)}=
R^{-1} \big( X_{i+1}(t) -R \big) R\,\ket{\Psi^{(i)}(t)}
\ef,
\ee
that is,
\be
\ket{\Psi^{(i+1)}(t)} = R\,\ket{\Psi^{(i)}(t)}
\label{eq.diheSER}
\ee
(more precisely, our reasoning only proves that 
$\ket{\Psi^{(i+1)}(t)} \propto R\,\ket{\Psi^{(i)}(t)}$, however,
as all the components are real-positive for $0<t<1$, and 
$\ket{\Psi^{(i+N)}(t)} \equiv \ket{\Psi^{(i)}(t)}$, the
proportionality factor must be 1).

Call $V$ the operator that applies a vertical symmetry to $\pi$, i.e.,
if $\{i,j\}$ is an arc of $\pi$, $\{N+1-j,N+1-i\}$ is an arc of
$V\pi$.  Clearly,
\begin{align}
V^2 &= \buno
\ef;
&
V R^{\pm 1} V &=
R^{\mp 1}
\ef;
&
V \beo_i V &=
\beo_{N-i}
\ef.
\end{align}
Then, as $X_i(t)^{-1} = X_i(t^{-1})$,
\be
\begin{split}
0 &= \big( X_i(t)-R \big)\,\ket{\Psi^{(i)}(t)}
= X_i(t) \big( R^{-1} - X_i(t^{-1}) \big) R\,\ket{\Psi^{(i)}(t)}
\\
&= -X_i(t) V \big( X_{N-i}(t^{-1}) - R \big) V R\,\ket{\Psi^{(i)}(t)}
\ef,
\end{split}
\ee
that is, $\ket{\Psi^{(N-i)}(t^{-1})} \propto VR\,\ket{\Psi^{(i)}(t)}$,
or, using (\ref{eq.diheSER}),
\be
\ket{\Psi^{(i)}(t^{-1})} \propto V\,\ket{\Psi^{(N+1-i)}(t)}
\ef.
\label{eq.diheSEV}
\ee
In fact the proportionality factor is $t^{-(n-1)}$ for the case
$\LP(2n)$, and $t^{-(N-1)}$ for the case $\LP^*(N)$, but this will only
be determined from the following Proposition \ref{prop.degT}.
Note that equation (\ref{eq.diheSEV}) corresponds to Claim 3
in~\cite{PdF04}. Also note that the covariance under vertical
reflection implies that our collection of solutions $\ket{\Psi^{(i)}(t)}$ of the
scattering equations (\ref{eq.qKZfund}) are related to the solutions
of the modified version in which the other direction of rotation is
considered, 
$X_i(t)\ket{\Psi}= R^{-1} \ket{\Psi}$.

At the light of the dihedral covariance discussed above, we can
concentrate on the case $i=1$ without loss of generality. All the
results that follow can be easily translated to arbitrary values of
$i$.  For short, call $\Psi(t)$ the vector $\Psi^{(1)}(t)$.

As $\beo_1^2=\beo_1$, the operators $\beo_1$ and $(\buno-\beo_1)$ are orthogonal
projectors, and a vector satisfies equation (\ref{eq.qKZfund}) if and
only if it satisfies the system of linear equations associated to
these two projectors
\begin{align}
\beo_1 \big( t\buno-R-(t-1)\beo_1 \big) \ket{\Psi(t)} 
&= 0
\ef;
\\
(1-\beo_1) \big( t\buno-R-(t-1)\beo_1 \big) \ket{\Psi(t)} 
&= 0
\ef.
\end{align}
Applying the rules of Temperley--Lieb Algebra leads to
\begin{align}
\beo_1 (\buno -R) 
\ket{\Psi(t)} 
&= 0
\ef;
\label{eq.452776A}
\\
(\buno -\beo_1) \big( t\buno-R \big) \ket{\Psi(t)} 
&= 0
\label{eq.452776}
\ef.
\end{align}
The first equation has the simplifying property of involving an
operator independent of $t$. The second equation is better
investigated in components
%
\be
\begin{split}
(1-\beo_1) \big( t\buno-R \big) \ket{\Psi(t)}
&=
(1-\beo_1) 
\Big(
\sum_{\pi} t \Psi(t;\pi) \ket{\pi}
-
\sum_{\pi} \Psi(t;\pi) \ket{R\pi}
\Big) 
\\
&=
(1-\beo_1) 
\Big(
\sum_{\pi} t \Psi(t;\pi) \ket{\pi}
-
\sum_{\pi} \Psi(t;R^{-1}\pi) \ket{\pi}
\Big) 
\\
&=
\sum_{\pi} 
\big( t \Psi(t;\pi)-
\Psi(t;R^{-1}\pi) \big)
(1-\beo_1) 
\ket{\pi}
\ef.
\end{split}
\ee
In this sum, a component $\pi$ with $1 \simX 2$ is annihilated by
$1-\beo_1$.  For (\ref{eq.452776}) to hold, as a component $\pi$ with
$1 \nsimX 2$ has coefficient $t \Psi(t;\pi)-\Psi(t;R^{-1}\pi)$, this
must be zero.
For later reference we record the previous observations as a
proposition.
\begin{proposition}
\label{prop.equiqKZ}
A necessary and sufficient condition for $\ket{\Psi(t)}$ to solve the
scattering equation at $i=1$, equation (\ref{eq.qKZfund}), is that
\begin{align}
\label{eq.qKZproj}
\beo_1 (\buno -R) 
\ket{\Psi(t)} 
&= 0
\ef;
\end{align}
and that, for any component $\pi$ such that $1 \nsimX 2$,
\begin{align}
\label{eq.PhiAnsatzRepeat}
t\,\Psi(t;\pi) 
&= \Psi(t;R^{-1} \pi)
\ef.
\end{align}
\end{proposition}
\noindent
This proposition will be employed in Section \ref{main-section}, in
the proof of Theorem \ref{main}. However, already at this point it is
useful for deducing a corollary of some importance.
\begin{corollary}
\label{coroll.degLB}
All the polynomial solutions $\Psi(t;\pi)$ of the scattering equation
have degree at least $n-1$ in the case $\LP(2n)$, and at least $N-1$
in the case $\LP^*(N)$.
\end{corollary}
\noindent
This comes from investigating the consequences on rainbow patterns of
equation (\ref{eq.PhiAnsatzRepeat}) in Proposition~\ref{prop.equiqKZ},

There are various possible derivations of the fact that a
normalisation of the solution exists such that the degree is
\emph{exactly} the one given above. In particular, one proof will
follow later on from our analysis (see Proposition \ref{prop.degT}).

\subsection{Further properties of the ground state}

\noindent
This section discusses some other properties of the ground states
$\Psi$ and $\Psi^{(i)}(t)$, that are not strictly necessary in the
derivation of the results of this paper, but are of interest in the
combinatorial connection with enumerations of FPLs, and allow to
understand some of the implications of our results.

As we stated above, it follows from the stochasticity of the
Hamiltonian and of the scattering matrix, the fact that the defining
equations involve integer coefficients, and from the Perron--Frobenius
theorem, that the components $\Psi(\pi)$ can be normalised to be
positive integers, and the components $\Psi^{(i)}(t;\pi)$ can be
normalised to be integer-valued polynomials in $t$, everywhere
positive in the open interval $t \in \;]0,1[$.

Much more than this is true. The components $\Psi(\pi)$ can be
normalised to `small' integers, in particular they remain all integers
if the entry of smallest component, corresponding to the rainbow
pattern, is set to 1.  Furthermore, the components $\Psi^{(i)}(t;\pi)$
can be normalised to be integer-valued polynomials in $t$, with all
positive coefficients, and in such a way that the coefficients are
still integers, when the entries of the various rainbow patterns (for
the various rotations) are set to monomials $t^j$. This means in
particular that the smallest normalisation that makes the
$\Psi(\pi)$'s all integers, and the smallest normalisation that makes
the $\Psi^{(i)}(t;\pi)$'s all integer-valued polynomials, are such
that $\Psi^{(i)}(t;\pi)|_{t=1} = \Psi(\pi)$, without proportionality
factors.  Finally, with this normalisation of the solution, the
maximum degree in the collections of polynomials coincides with the
lower bounds given in Corollary~\ref{coroll.degLB}.

All of the claims of the previous paragraph follow from the
Razumov--Stroganov correspondence, either in its original form (for
the claims concerning $\Psi(\pi)$) or in the refined version discussed
here (for the claims concerning $\Psi^{(i)}(t;\pi)$), through the
interpretation of the ground state components as enumerations of
FPL's. Some of the claims have alternative, purely
algebraic proofs, which are mostly based on the observation that the
vectors $\ket{\Psi^{(i)}(t)}$ are a special case of
$\ket{\Psi(z_1,\dots, z_N)}$, the ground state of the fully
inhomegenous $O(1)$ Dense Loop Model.
More precisely,
\be
\label{eq.reltoqKZ}
\ket{\Psi^{(i)}(t)} \propto
\ket{\Psi
\big(1,\ldots,1, z_{i}=\smfrac{qt+q^{-1}}{q+q^{-1}t}, 1, \ldots,1
\big) }
\ef.
\ee
The claim on the degree in $t$ can be deduced by the analogous claim
on the degree in the $z_i$'s, provided in \cite{pdf-pzj1} for the
$\LP(2n)$ case, and in \cite{pdf-pzj-jbz06} for the $\LP^*(N)$
case. As a byproduct, one also obtains the fact that rainbow pattern
components are of the form $t^j$.
Furthermore, also the claim on the integrality of the components
$\Psi(\pi)$ can be proven (see e.g.\ \cite[sec.\;4.4]{PZJthesis}).

We have seen that, for $t=1$, the ground states of the scattering
matrix, $\Psi^{(i)}(t;\pi)$, reduce to the one of the Hamiltonian.
Another interesting specialization is $t=0$ (or its symmetric $t \to
\infty$).  From equation (\ref{eq.PhiAnsatzRepeat}) in
Proposition~\ref{prop.equiqKZ}, in this limit $\ket{\Psi_N^{(i)}(0)}$
lies in the image of $\beo_{i-1}$, i.e., for $\pi \in \LP(N)$ (or
$\LP^*(N)$), $\Psi_N^{(i)}(0,\pi) = 0$ if $i-1 \nsimX i$ in
$\pi$. This subset of link patterns is in natural correspondence with
the full set $\LP(N-2)$ (or $\LP^*(N-2)$), upon removing the short arc
between $i-1$ and $i$.  As a consequence of a recursion relation
satisfied by $\ket{\Psi(z_1,\dots, z_N)}$ and proven in
\cite{pdf-pzj1} (Theorem 3), this specialisation and correspondence
produces a vector that satisfies equation (\ref{RS-eq0}) at size
$N-2$. The solution of this equation with non-negative entries is
unique, and the customary fixing of the normalisation, e.g.\ for the
rainbow pattern containing the arc $\{i-1,i\}$, allows to fix the
proportionality factor to~1:
\be
\label{recursion-paul-philippe}
\Psi_N^{(i)}
\Big( 0,
\setlength{\unitlength}{6pt}
\begin{picture}(8,2.8)
\put(0,0){\includegraphics[scale=1.2]{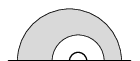}}
\put(2.2,-1){${}_{i-1}\; {}_i$}
\put(1.8,0.8){$\pi$}
\end{picture}
\Big)
=
\Psi_{N-2}^{(i)}(1,\pi)
\ef.
\ee

\section{Fully-packed loops}
\label{sec.FPLdefs}

\noindent
In this section we present the relevant facts on Fully-Packed Loops
that we will need in the proofs of our results. Section
\ref{ssec.fplBasic} presents the main definitions. Section
\ref{ssec.multic} introduces a family of domains, that we call
\emph{dihedral domains}, for which our theorems do hold.  These
domains have a common behaviour under \emph{gyration}, a bijection
introduced by Wieland \cite{wie} for FPL's in the square. This
operation is analysed in detail in Section \ref{sec.WieRemind}.



\subsection{Fully-packed loops: basic definitions}
\label{ssec.fplBasic}

\noindent
Consider a graph $\Lambda$ in which all the vertices have degree in
the set $\{1,2,4\}$. If the degree is 2 or 4, we say that the vertex
is \emph{internal}, while if the degree is 1 we say that the vertex is
\emph{external}, and that the incident edge is an \emph{external edge}.
Obvious parity reasons force to have an even number of external edges.

We say that $\Lambda$ is \emph{outer-planar} if it admits a planar
embedding into a disk, with all the external vertices on the boundary
of the disk. Such an embedding defines univocally the \emph{faces} of
the domain, and induces a cyclic ordering of the external vertices and
edges.

We will be interested in maps $\phi: E(\Lambda) \to \{b,w\}$,
naturally identified with bicolourations of the edge-set of the graph
(and $b$, $w$ stand for black, white). The map $\phi$ is a
\emph{fully-packed loop} (FPL) configuration if all the internal
vertices of degree 4 have two black and two white adjacent edges, and
all the vertices of degree 2 have one black and one white adjacent
edge.  The map $\phi$ has \emph{boundary condition} $\tau$ if its
restriction to the external edges is $\tau$. We define
$\fpl(\Lambda;\tau)$ the set of FPL's on $\Lambda$, with boundary
condition $\tau$. We call $\overline{\tau}$ the boundary condition
with black and white exchanged. For $\phi \in \fpl(\Lambda;\tau)$, the
\emph{conjugate} configuration $\sigma(\phi) \in
\fpl(\Lambda;\overline{\tau})$ is the one obtained interchanging black
and white.  If $\Lambda$ is outer-planar, there are two specially
important boundary conditions, $\tau_+ = (b,w,b,w,\ldots,b,w)$ and
$\tau_- = \overline{\tau_+} = (w,b,w,b,\ldots,w,b)$, that we call
\emph{alternating boundary conditions}. We also use the shortcuts
$\fpl_+(\Lambda) = \fpl(\Lambda;\tau_+)$ and
$\fpl_-(\Lambda) = \fpl(\Lambda;\tau_-)$.

The black and the white edges of $\phi$, separately, form a collection
of cycles and of open paths with endpoints at vertices of degree $1$
or $2$.
The open paths of each given colouration connect external vertices and
internal vertices of degree $2$ among themselves, thus inducing two
pairings, $\pi_b$ and $\pi_w$ for black and white respectively, and
two integers, $\ell_b$ and $\ell_w$, for the number of black and white
cycles in $\phi$. Call $\ell=\ell_b+\ell_w$.

If $\Lambda$ is outer-planar, and has no vertices of degree 2, and
$\tau$ has $2n_b$ black and $2n_w$ white entries, $\pi_b$ and $\pi_w$
are link patterns in $\LP(2n_b)$ and $\LP(2n_w)$, i.e.\ they have no
crossings. If $\Lambda$ has a marked face, $\pi_b$ and $\pi_w$ are
punctured link patterns in $\LP^*(2n_b)$ and in $\LP^*(2n_w)$, with
the puncture located inside the marked face. If it has exactly one
vertex of degree 2, and $\tau$ has $2n_b-1$ black and $2n_w-1$ white
entries, $\pi_b$ and $\pi_w$ are punctured link patterns in
$\LP^*(2n_b-1)$ and $\LP^*(2n_w-1)$, with the puncture located at this
vertex.

The case corresponding to the original Razumov--Stroganov
correspondence is the one in which $\Lambda$ is a square $n \times n$
grid, with $n^2$ internal vertices of degree 4, and $4n$ external
vertices, adjacent to the internal vertices in the obvious way. Note
that, for $n$ sufficiently large, the outer-planar embedding is
unique, up to reflection symmetry, and to permutation of the pairs of
external edges adjacent to corners.


\subsection{Dihedral domains: definition and properties} 
\label{ssec.multic}

\noindent
A tool in the combinatorial investigation of FPL's, started with the
work of Wieland \cite{wie}, is that on the $n \times n$ square with
alternating boundary conditions
there exist two bijections $H_{\pm}$
between $\fpl_+(\Lambda)$ and $\fpl_-(\Lambda)$, that preserve the
link pattern up to a rotation, which
imply relations between the enumerations of FPL's with a given
triple $(\pi_b,\pi_w,\ell)$. 

This property of the maps $H_\pm$ has been crucial for the proof of
the Razumov--Stroganov correspondence in \cite{usRS},
where it was mentioned that analogues of
the maps $H_\pm$ exist on a more general class of outer-planar domains
characterized by the following two properties:
\begin{itemize}
\item the internal faces have at most $4$ sides,
\item the faces formed by joining two consecutives
  external edges have at most $3$ sides.
\end{itemize}
In the present paper
we exhibit a family of graphs, that we call \emph{dihedral domains},
that satisfy these properties. All our results holds on this class of
domains.


In \cite{usprep}
we shall perform a classification of \emph{all} the domains presenting
a version of the Wieland Theorem, without any \emph{a priori}
assumption (for example, without assuming outer-planarity).

We produce two families of dihedral domains. Domains of the first kind
present the version of the Razumov--Stroganov correspondence on
$\LP(2n)$ and, in certain circumstances, also the more refined version
on $\LP^*(2n)$.  Domains of the second kind present the version of the
correspondence on $\LP^*(2n-1)$.

\begin{figure}[!tb]
\[
\includegraphics[scale=2]{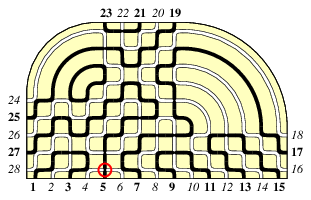}
\qquad
\includegraphics[scale=2]{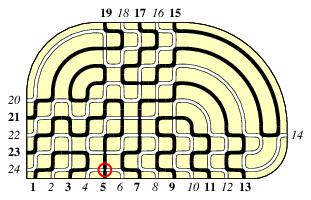}
\]
\caption{\label{fig.multicex}Typical domains $\Lambda$, with
  two corners (left) and one corner (right),
  together with a configuration $\phi \in \fpl_+(\Lambda)$.  The
  domain on the left is $\Lambda(15,9;0,0,6,4)$, and has $2n=14$ and
  $L=15$.  The domain on the right is $\Lambda(15,9;0,2,6,4)$, and has
  $2n=12$ and $L=24$. In both cases, $h(\phi) = 5$, because on the
  reference line the position of the only $c$-type tile is adjacent to
  the 5-th external edge.}
\end{figure}


A domain of the first
kind
has six integer parameters, $\Lambda =
\Lambda(L_x,L_y;a_1,a_2,a_3,a_4)$.  Take a rectangular $L_x \times
L_y$ portion of the square grid, and call $A_1, \ldots, A_4$ the four
corners, in counter-clockwise order starting from the bottom-left one.
From the corner $A_i$, cut away an $a_i\times a_i$ square portion of
the grid. The parameters must be such that the removed squares do not
overlap, although they might share a portion of the perimeter (e.g.,
$a_1+a_2 = L_x$ is admitted, but $a_1+a_2 = L_x+1$ is not).  Each
removed square will cut $2a_i$ internal edges, which are connected
pairwise, starting from the resulting concave angle. The number of
external edges is $2(L_x+L_y-a_1-a_2-a_3-a_4)$, and the parameter $2n$
of the set of link patterns will be given by
\be
2n=L_x+L_y-a_1-a_2-a_3-a_4 \ef.
\ee
Domains defined by this procedure have all faces with at most $4$
sides.  They are constituted of a finite number of bundles of parallel
lines (at least 2 and at most 6 bundles), that cross each other
forming portions of the square grid. These lines are the geometric
structures to which, in the general framework of the integrable
6-Vertex Model, are associated spectral parameters. One line is a
\emph{boundary line} if it is adjacent to some external edges.

Define a \emph{corner} as a vertex adjacent to two external
edges, and call $c_\Lambda$ the number of corners in $\Lambda$, which
also coincides with the number of boundary lines. For a
face with $\ell$ sides, the \emph{curvature} of the face is
$4-\ell$. The total curvature of the domain, $d_\Lambda$, is the sum
of the curvatures of all the faces.  The number of corners
$c_\Lambda$
is exactly $4$ minus the number of non-zero $a_i$'s, and is related to the
total curvature by an Euler formula, $c_\Lambda=4-d_\Lambda$.
Examples of dihedral domains with $c_\Lambda=2$ and $c_\Lambda=1$ are
shown in Figure \ref{fig.multicex}, together with examples of FPL.

\begin{figure}
\[
\includegraphics[scale=2]{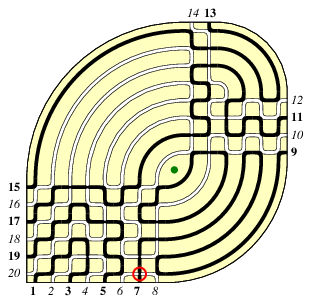}
\qquad
\includegraphics[scale=2]{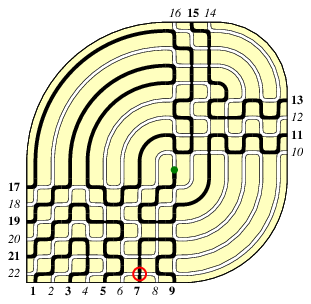}
\]
\caption{\label{fig.multicexPunct} Left: a typical domain $\Lambda$ of
  the first kind, with one face of two sides, showing
  Razumov--Stroganov correspondence for punctured link patterns
  $\LP^*(2n)$. This domain is $\Lambda(14,14;0,6,4,8)$, and has
  $2n=10$ and $L=20$. It has one face with two sides because
  $L_x=L_y=a_2+a_4$.  Right: a typical domain $\Lambda$ of the second
  kind, with one vertex of degree 2, showing Razumov--Stroganov
  correspondence for punctured link patterns $\LP^*(2n-1)$. This
  domain is $\Lambda(15,14;0,6,4,8)$, and has $2n-1=11$ and $L=22$.
  It allows for a vertex of degree 2, because $L_x-1=L_y=a_2+a_4$.
  The puncture is marked in green.  Typical configurations $\phi \in
  \fpl_+(\Lambda)$ are shown.  both with $h(\phi) = 7$.}
\end{figure}

As domains with $c_\Lambda=0$ have refined enumerations of FPL's which
are easily related to the ones of a reduced domain, in the following
we restrict to the `irreducible' case in which there is at least one
corner (thus $c_\Lambda\geq 1$, or equivalently $d_\Lambda \leq 3$).

Because of the rotational invariance of the parameter space, we can
set $a_1=0$ without loss of generality, and set the corner $A_1$ as
the \emph{reference corner} of the domain. The external edges are
labeled in counter-clockwise order, starting from this corner, and the
\emph{reference side} of the domain is the boundary line adjacent to
the external edges with labels $1,2,\ldots,$ up to when another corner
is reached. We will call $L$ the length of the reference side,
i.e.\ the number of vertices along this side. The dependence of $L$
from our parameters depends on which corners are present, as $L=L_x$
if $a_2=0$, $L=L_x+L_y-2a_2$ if $a_2>0$ and $a_3=0$, and so on.

Under the restriction to $c_\Lambda\geq 1$, our domains have at most
one face with $1$ or $2$ sides.  The set of domains with one such face
coincides with the set of domains that allow for the $\LP^*(2n)$
version of the Razumov--Stroganov correspondence, and the puncture is
located inside this special face.  A dihedral domain of this type,
together with a FPL configuration, is shown in Figure
\ref{fig.multicexPunct}, left.

Dihedral domains of the second kind are a variant of a subclass of
domains of the first kind, that have an internal edge $e$ for which
both adjacent faces have at most 3 sides. This edge is splitted in
two, and a vertex of degree $2$ is produced. Note that, after
this procedure, still the adjacent faces have at most 4 sides, as
required. Note also that, again by reasonings of curvature, this
procedure cannot be applied more than once (indeed, in the Euler
formula, a vertex of degree $4-d$ carries $d$ units of curvature). The
parameter $2n-1$ of the set of link patterns $\LP^*(2n-1)$ is given by
\be
2n-1=L_x+L_y-a_1-a_2-a_3-a_4
\ef,
\ee
and the puncture of the link patterns is located at the vertex of
degree 2. A dihedral domain of second kind, together with a FPL
configuration, is shown in Figure \ref{fig.multicexPunct}, right.


Given a dihedral domain $\Lambda$, we will be interested in the set
$\fpl(\Lambda)$ of FPL configurations with alternating boundary
conditions.  Depending on the parity of the labels of the external
black edges, the set $\fpl(\Lambda)$ splits in two equinumerous
subsets, $\fpl_+(\Lambda)$ and $\fpl_-(\Lambda)$, which contains FPL
such that the external edge with label 1 is black or white. The
operation $\sigma$, which consists in inverting the colouration of all
the edges, is an involution in $\fpl(\Lambda)$, which provides a
trivial bijection between $\fpl_+(\Lambda)$ and $\fpl_-(\Lambda)$.

We recall that at every vertex of degree 4, given an ordering of the
two crossing lines, there are 6 possible local configurations (this is
the origin of the name ``6-Vertex Model''), pairwise related by the
involution $\sigma$
\be
\label{eq.defabc}
\setlength{\unitlength}{20pt}
\raisebox{-4mm}{
\begin{picture}(13.5,1.4)(0,-0.2)
\put(0,0){\includegraphics[scale=4]{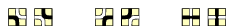}}
\put(0,0.4){$a:$}
\put(5,0.4){$b:$}
\put(10,0.4){$c:$}
\end{picture}
}
\ee
We will say that a FPL $\phi$ has a tile of type $a$, $b$ or $c$ at a
given vertex of degree-4, when the ordering of the two lines is
unambiguous. This is in particular the case for vertices adjacent to
an external line.

It is easy to see that, because of the alternating boundary conditions, each $\phi
\in \fpl(\Lambda)$ is forced to have a single $c$-tile on any of the
boundary sides, with all $a$-tiles and $b$-tiles along the line, at
the right and the left of the $c$-tile, respectively.
In particular, there is a single $c$-tile at position $1 \leq h(\phi)
\leq L$ along the reference side. We call $h(\phi)$ the
\emph{refinement position} of the FPL $\phi$ (it is understood,
w.r.t.\ the reference side).

We call $\fpl_b(\Lambda)$ and $\fpl_w(\Lambda)$ the subsets of
$\fpl(\Lambda)$ consisting of FPL such that the edge incident to the
refinement position is black or white.  The operation $\sigma$ is also
trivially a bijection between $\fpl_b(\Lambda)$ and $\fpl_w(\Lambda)$.

For a dihedral domain $\Lambda$, we will denote with $\LP(\Lambda)$ 
the set $\LP(2n)$, $\LP^*(2n)$ or $\LP^*(2n-1)$, depending from the
kind of domain, and for the appropriate value of $n$. The only
ambiguity holds for those domains of the first kind
for which both the $\LP(2n)$ and the $\LP^*(2n)$ correspondences are
viable. As the latter is a stronger version of the former, we can
assume without loss of generality that $\LP(\Lambda)=\LP^*(2n)$ in
such a case.

\begin{figure}[!tb]
\setlength{\unitlength}{10pt}
\begin{picture}(38,41.5)
\put(0,21){\includegraphics[scale=2]{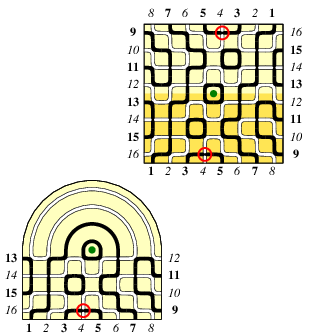}}
\put(17,21){\includegraphics[scale=2]{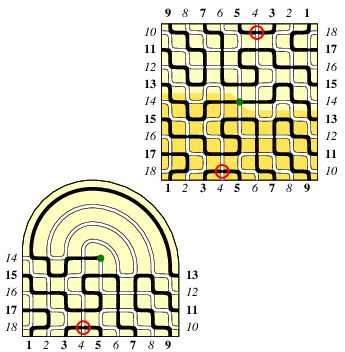}}
\put(0,0){\includegraphics[scale=2]{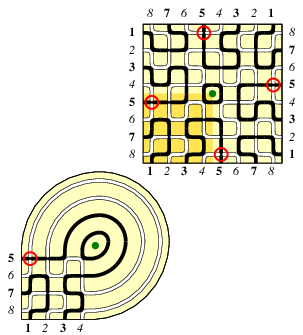}}
\put(17,0){\includegraphics[scale=2]{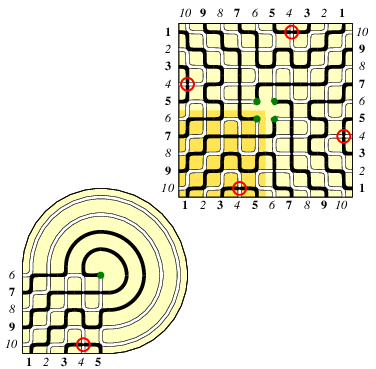}}
\put(0,31){\includegraphics[scale=0.66]{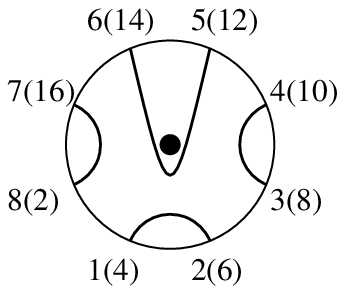}}
\put(29.7,23.7){\includegraphics[scale=0.66]{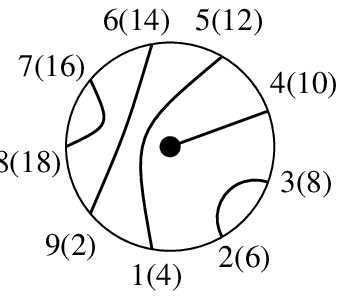}}
\put(0,10){\includegraphics[scale=0.66]{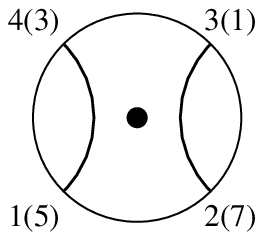}}
\put(29.7,2.7){\includegraphics[scale=0.66]{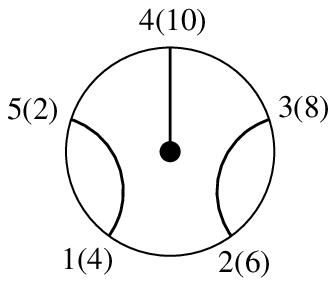}}
\end{picture}
\caption{\label{fig.symclass}Dihedral domains associated to symmetry
  classes of ASM that present dihedral Razumov--Stroganov
  correspondence. Top-left and top-right: HTASM of even and odd size,
  respectively; bottom-left: QTASM of size $4n$; bottom-right:
  quasi-QTASM of size $4n+2$. In the link patterns, the notation
  $i\,(j)$ stands for the fact that the link pattern termination is
  the image on the FPL of an external edge, which is the $i$-th one
  counting from the refinement position, and of the colour of the
  refinement position, and the $j$-th one (including both colourings)
  counting from the reference corner.}
\end{figure}

We conclude the present section by observing that some particular
dihedral domains correspond to symmetry classes of FPL's (and thus of
Alternating Sign Matrices) that have been already studied in the
literature (see \cite{oneroof}, and \cite{duchon} for quasi-QTASM).
Figure \ref{fig.symclass} gives illustrations of these classes.

\begin{description}
\item[ASM:] When all $a_i$ are zero, and $L_x=L_y$, the domain is the
  square, and the FPLs are in bijection with all Alternating Sign
  Matrices. Note that conversely, when all $a_i$ are zero and $L_x
  \neq L_y$, $\fpl(\Lambda)=\emptyset$.
\item[HTASM:] For any two of the $a_i$'s equal to 0, and the other two equal
  to $n$, and $L_x=L_y=2n$, we get Half-Turn Symmetric ASM (HTASM) of
  size $2n$, while for $a_1=a_3=0$, $a_2=a_4=n$ and $L_x-1=L_y=2n$, as
  we have an edge adjacent to two triangular faces, we can construct a
  dihedral domain of the second kind, which corresponds to HTASM of
  size $2n+1$ (see Figure \ref{fig.symclass}, top).
\item[QTASM:] For $a_1=0$, $a_2=a_3=a_4=2n$ and $L_x=L_y=4n$ we get
  Quarter-Turn Symmetric ASM (QTASM) of size $4n$, while for
  $L_x-1=L_y=2n$ we can construct a dihedral domain of the second
  kind, which leads to the so called quasi-QTASM of size $4n+2$ (see
  Figure \ref{fig.symclass}, bottom). Note that we do \emph{not}
  present dihedral domains corresponding to QTASM of size $4n\pm 1$.
\end{description}

\subsection{Wieland gyration: a reminder of facts}
\label{sec.WieRemind}

\noindent
For $\ell \geq 1$, call $C_{\ell}$ the \emph{cycle graph} with $\ell$
vertices, i.e.\ the graph with vertices $\{1,2,\ldots,\ell\}$ and
edges $\{(1,2),(2,3),\ldots,(\ell-1,\ell),(\ell,1)\}$. Consider
bicolourations $\phi$, in black and white, of a cycle $\gamma \cong
C_{\ell}$, for $\ell \leq 4$.  A cycle with $\ell \leq 2$ may also be
`punctured'.

\begin{figure}
\[
\includegraphics[scale=1]{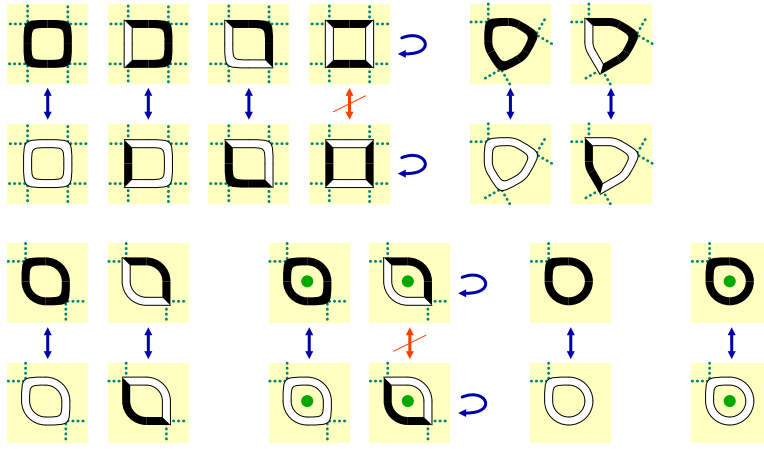}
\]
\caption{\label{fig.wie1234}Illustration of the action of local
  gyration $H_\gamma$, on all the possible colourings of cycles with
  length $\ell \leq 4$. Punctured cycles have a green mark.}
\end{figure}

We define the \emph{local gyration} $H_\gamma$ as the involutive map
on these colourings, that inverts the colouration of the edges, except
when $\ell=4$, or $\ell=2$ and the cycle is punctured, and the edges
are coloured in an alternate way, for which the map is the identity
(see Figure \ref{fig.wie1234}).
This involution has a number of important properties. It preserves the
\emph{connectivity of open paths} of a given colour, i.e., if $v_1$
and $v_2$ are the endpoints in $\gamma$ of a black (or a white) open
path w.r.t.\ $\phi$, they are still the endpoints of a black (or a
white) open path w.r.t.\ $H_{\gamma} \phi$. Furthermore, the two paths
are \emph{homotopically equivalent}, i.e.\ they can be deformed
continuously the one into the other, throughout the region inside the
cycle (and without crossing the puncture, if present).
Then, if and only if $\phi$ consists of a monochromatic cycle, also
$H_{\gamma} \phi$ has this property (and
both $\phi$ and $H_{\gamma} \phi$ consist of cycles than encircle the
puncture, if present).  Finally, the degree of a given colour at each vertex is
\emph{complemented}, i.e.\ $\deg_b(v,\phi) = 2 - \deg_b(v,H_{\gamma}
\phi)$ and $\deg_w(v,\phi) = 2 - \deg_w(v,H_{\gamma} \phi)$.

This last property implies that, for FPL configurations $\phi$ on a
domain $\Lambda$, if $H_{\gamma}$ is performed on a cycle $\gamma
\subseteq \Lambda$, the degree constraints at the internal vertices
are not necessarily preserved. However, suppose that the domain
$\Lambda$ has no external vertices, i.e.\ all of its vertices have
degree $2$ or $4$.
Assume that $\Lambda$ admits a \emph{cycle decomposition}
$\Gamma=(\gamma_1,\dots,\gamma_m)$, i.e.\ a partition of the edge-set
into graphs $\gamma_i$ isomorphic to cycle graphs with lengths
$|\gamma_i| \leq 4$. Then, because of the complementation property,
the involution $H_{\Gamma} = \prod_i H_{\gamma_i}$ preserves the
degree constraints at all the vertices, i.e.\ sends FPL into FPL (note
that the product is unambiguous because $H_{\gamma}$ and $H_{\gamma'}$
commute if $\gamma$ and $\gamma'$ are edge-disjoint cycles).

These observations imply the following
\begin{lemma}
\label{connect-preserv}
Let $\Lambda$ be a domain with no external vertices, and $2n$ internal
vertices of degree 2. Let $\Gamma=(\gamma_1,\dots,\gamma_m)$ be a
cycle decomposition of $\Lambda$, and let $\phi \in \fpl(\Lambda)$ and
$\phi' = H_{\Gamma}(\phi)$.  The following properties hold:
\begin{itemize}
\item[a)] If two degree-2 vertices $v_1$, $v_2$ are endpoints of a
  black path of $\phi$, then they are endpoints of a black path of
  $\phi'$. The same holds for white paths. Thus, if $\phi$ has
  $(\pi_b, \pi_w)$ black and white pairings of the endpoints, $\phi'$
  has the same black and white pairings $(\pi_b, \pi_w)$. Furthermore,
  if $\phi$ has $\ell = \ell_b+\ell_w$ cycles, also $\phi'$ has $\ell$
  cycles.

\item[b)] For each $e \in \gamma_i$, unpunctured and with $|\gamma_i|
  \leq 3$, if $\phi(e)=b$ then $\phi'(e)=w$, and vice versa.

\item[c)] If $\Lambda$ has a planar embedding which is outer-planar
  w.r.t.\ all the degree-2 vertices, then $\pi_b$, $\pi_w$ are both
  link patterns in $\LP(2n)$ (i.e., the pairings of a given colour are
  non-crossing w.r.t.\ the outer-planar embedding). The link patterns
  $\pi_b$, $\pi_w$ are preserved by $H_{\Gamma}$ also as punctured
  link patterns, in $\LP^*(2n)$, if the puncture is outside the cycles
  of $\Gamma$,
  or it is inside a face $\gamma_i \in \Gamma$ with $|\gamma_i| = 1$ or
  $2$. If $\ell^* = \ell_b^*+\ell_w^*$ is the number of cycles
  encircling such a puncture, this number is also preserved by
  $H_{\Gamma}$.

\item[d)] If $\Lambda$ has a planar embedding which is outer-planar
  w.r.t.\ all the degree-2 vertices except one, then $\pi_b$, $\pi_w$
  are link patterns in $\LP^*(2n-1)$, with the puncture being the only
  degree-2 vertex not on the boundary, and are preserved by
  $H_{\Gamma}$ also as punctured link patterns. In such a case,
  $\ell^* =0$ both in $\phi$ and in $\phi'$.
\end{itemize}
\end{lemma}

\proof
The proof consists just in translating all the established local
properties of gyrations $H_{\gamma}$ at the global level of
monochromatic cycles and open paths of $\Lambda$, observing that a
cycle on $\Lambda$ is either one of the $\gamma_i \in \Gamma$, or is
the concatenation of two or more open paths contained in the
$\gamma_i$'s, and that an open path is the concatenation of one or
more open paths contained in the $\gamma_i$'s.
\qed

\bigskip
\noindent
The class of graphs $\Lambda$ allowing for a cycle decomposition, and
thus a gyration operation $H_{\Gamma}$ with the properties above, is
very large.  Furthermore, for outer-planar graphs, in principle
multiple punctures can be introduced simultaneously (e.g., at several
inequivalent points which are all outside the cycles $\gamma_i$).
However, in such a case the properties deduced from Lemma
\ref{connect-preserv} are not specially useful.

On the contrary, if we require the presence of \emph{two inequivalent}
operations, the deduced properties become much more interesting, but
the class of domains becomes much more narrow, and no more than one
puncture can be introduced.  The dihedral domains described in the
previous section have the characteristics above.

\begin{figure}
\[
\setlength{\unitlength}{10pt}
\begin{picture}(29,9)
\put(0,0){\includegraphics[scale=2]{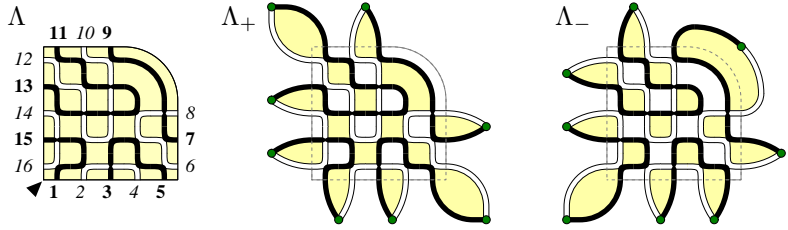}}
\put(0.25,7.75){$\Lambda$}
\put(8.25,7.75){$\Lambda_+$}
\put(20.75,7.75){$\Lambda_-$}
\end{picture}
\]
\caption{\label{fig.LLpLm}Left: a typical domain $\Lambda$.  If the
  reference corner is at the bottom left, the reference side is the
  bottom row, and the FPL configuration shown is in $\fpl_+(\Lambda)$.
  Middle and right: the construction of $\Lambda_+$ and $\Lambda_-$,
  respectively.}
\end{figure}

If $\Lambda$ is a graph with all vertices of degree 1, 2 or 4,
equipped with a planar embedding, outer-planar w.r.t.\ the vertices of
degree 1, we have two associated graphs $\Lambda_{\pm}$ with no
external vertices, for the two possible pairing of consecutive
external vertices (see Figure \ref{fig.LLpLm} for an example). A FPL
$\phi$ on $\Lambda$ induces also FPL's
both on $\Lambda_{+}$ and $\Lambda_{-}$ if and only if we have
alternating boundary conditions on $\Lambda$.

Consider the case in which $\Lambda$ is a dihedral domain with $2N$
external vertices.
As explained in Section \ref{ssec.multic}, we have a reference corner,
a reference side, associated to an external line, and a natural
counter-clockwise labeling of the external edges.

The graphs $\Lambda_{\pm}$ inherit the planar embedding from
$\Lambda$. Moreover, through the induced structure of two-dimensional
cell complex, they have obvious cycle decompositions $\Gamma_{\pm}$.
The two gyration operations, $H_{\Gamma_+}$ on $\Lambda_+$ and
$H_{\Gamma_-}$ on $\Lambda_-$, will be just called $H_+$ and $H_-$ in
the following.

All the external edges $e$ of $\Lambda$ are in $\Lambda_{\pm}$ within
cycles $\gamma_i$ of length at most 3, so that the property {\it b} of
Lemma \ref{connect-preserv} applies. This is important in order to
establish that the operations $H_{\pm}$, which act on $\Lambda_{\pm}$,
when reinterpreted at the level of $\Lambda$ induce two involutions
on $\fpl(\Lambda)$, which are bijections between $\fpl_+(\Lambda)$ and
$\fpl_-(\Lambda)$.
Another important aspect is that the property of preserving $\pi_b$,
$\pi_w$, $\ell$ and (if applicable) $\ell^*$ on $\Lambda_{\pm}$ has a
counterpart also on~$\Lambda$.


In order to see this, at the level of generality required for the
refined Razumov--Stroganov correspondence, we introduce a general
class of maps that associate a link pattern to a FPL configuration.
For $v$ an external vertex of $\Lambda$ such that the adjacent edge
has boundary condition black, we let $\Pi(\phi,v)$ the map that
associates to $\phi$ the black link pattern, with label $1$ associated
to vertex $v$, and so on up to $N$ with increasing labels in
counter-clockwise order.  Two special cases of $\Pi(\phi,v)$ are
mostly used in the following. For $\phi \in \fpl_+(\Lambda)$, we
define $\Pi_+(\phi) = \Pi(\phi,1)$, i.e.\ the external vertex with
label $1$ is the first one of the reference side, next to the
reference corner.  This is the `ordinary' map used in the
correspondence, i.e.\ $\Pi_+(\phi) = \pi(\phi)$ in the notations of
the introduction.  Then, for $\phi\in \fpl_b(\Lambda)$, we set
$\Pi_b(\phi) = \Pi(\phi,h(\phi))$, i.e.\ the counting starts from the
refinement position. This is the `new' map introduced in this paper,
i.e.\ $\Pi_b(\phi) = \tilde{\pi}(\phi)$ in the notations of the
introduction (see Figure \ref{fig.expitilde}).

Remark that, as we have alternating boundary conditions, if the $v$-th
external edge is black, the set of black external edges is
$\{\ldots,v-2,v,v+2,\ldots\}$, and
as a consequence
\be
\label{eq.rotaPiV}
\Pi(\phi,v+2) = R^{-1} \Pi(\phi,v)
\ef.
\ee
Both maps $H_{\pm}$ are involutions. They do not commute among
themselves, and clearly $(H_+ H_-)^{-1} = H_- H_+$, so the only
irreducible elements in the generated monoid are strings of the form
$\cdots H_- H_+ H_- H_+ \cdots$.  It is thus convenient
to call $H: \fpl(\Lambda)\rightarrow \fpl(\Lambda)$ the map
\be
\label{H-map-def}
H(\phi)
=
\left\{
\begin{array}{ll}
H_+(\phi) & \textrm{if $\phi\in\fpl_+(\Lambda)$;}\\
H_-(\phi) & \textrm{if $\phi\in\fpl_-(\Lambda)$.}
\end{array}
\right.
\ee
Acting on $\fpl_+(\Lambda)$ we have
\begin{align}
H^k(\phi) 
&= 
\underbrace{H_{(-1)^{k-1}} \cdots H_- H_+}_{k} (\phi)
\ef;
&
H^{-k}(\phi) 
&= 
\underbrace{H_{(-1)^{k}} \cdots H_+ H_-}_{k} (\phi)
\ef;
\end{align}
and the irreducible elements of the monoid are just the powers 
$\{ H^k \}_{k \in \mathbb{Z}}$.

The following lemma, due to Wieland in the case of the square
\cite{wie}, is deduced easily from Lemma \ref{connect-preserv}
and the construction of $\Lambda_{\pm}$ from a domain $\Lambda$
described above, and will be at basis of the reasonings in Section
\ref{main-section}  
\begin{lemma}[Wieland half-gyration lemma]
Let $\Lambda$ be a dihedral domain, $\phi \in \fpl(\Lambda)$ and $v$ an
external vertex of $\Lambda$ adjacent to a black edge of $\phi$. Then
\be
\Pi(H(\phi),v+1) = \Pi(\phi,v)
\ef,
\label{wiel-gyr-fplA}
\ee
and in particular
\be
\label{wiel-gyr-fpl}
\Pi(H^2(\phi),v)= R\,\Pi(\phi,v)
\ef.
\ee
\end{lemma}
\proof Say that $v$ is odd, so that $\phi\in\fpl_+(\Lambda)$ and $H
\equiv H_+$. In $\Lambda_+$, the vertex $v$ is paired to $v+1$. After
the application of $H$, the link pattern $\pi_b$ is preserved on
$\Lambda_+$. We can split back the vertices, to obtain a configuration
in $\fpl_-(\Lambda)$. In particular, as we can apply the property 
{\it b} of Lemma \ref{connect-preserv} to the edges adjacent to $v$
and $v+1$, in $H(\phi)$ the black edge is now adjacent to $v+1$, from
which equation (\ref{wiel-gyr-fplA}) follows. If $v$ is even, the
reasoning is analogous. In particular, because of the staggered choice
of equation (\ref{H-map-def}), it is still true that $v$ is paired to
$v+1$.  Equation (\ref{wiel-gyr-fpl}) comes from applying
(\ref{wiel-gyr-fplA}) twice and then (\ref{eq.rotaPiV}):
$\Pi(\phi,v)=\Pi(H^2(\phi),v+2)=R^{-1}\Pi(H^2(\phi),v)$.
\qed

\bigskip \noindent 
A stronger version of the lemma (that we do not use in this paper),
still due to Wieland in the case of the square, considers the triple
$(\pi_b,\pi_w,\ell)$, or, in case of $\LP(\Lambda)=\LP^*(2n)$, the
quadruple $(\pi_b,\pi_w,\ell,\ell^*)$. Let
${\bm \Pi}(\phi,v,v')= 
\xi^{\ell(\phi)} \eta^{\ell^*(\phi)}
\;
\Pi(\phi,v) \otimes \Pi(\sigma \phi,v')$
(with $\ell^*(\phi)\equiv 0$ if not applicable). 
Then we have
\begin{lemma}
Let $\Lambda$ be a dihedral domain, $\phi \in \fpl(\Lambda)$ and $v$,
$v'$ external vertices of $\Lambda$ adjacent to a black and a white
edge of $\phi$, respectively. Then
\be
{\bm \Pi}(H(\phi),v+1,v'-1) = {\bm \Pi}(\phi,v,v')
\ef,
\ee
and in particular
\be
\label{wiel-gyr-fpl.3}
{\bm \Pi}(H^2(\phi),v,v') =
(R \otimes R^{-1}) \circ
{\bm \Pi}(\phi,v,v')
\ef.
\ee
\end{lemma}

\noindent
The Wieland gyration theorem of \cite{wie}, generalised to dihedral
domains, then follows as a corollary of equations (\ref{wiel-gyr-fpl})
and (\ref{wiel-gyr-fpl.3}).
\begin{theorem}[Wieland gyration theorem]
Let $\Lambda$ be a dihedral domain, and $\Psi_\Lambda(\pi)$
the number of FPL $\phi \in \fpl_+(\Lambda)$ such that 
$\Pi_+(\phi) = \pi$.
Then
\be
\Psi_\Lambda(\pi) =
\Psi_\Lambda(R^{-1} \pi)
\ef.
\ee
Similarly, for $\Psi_\Lambda(\pi_b,\pi_w,\ell,\ell^*)$ the number of
FPL $\phi \in \fpl_+(\Lambda)$ with the corresponding quadruple
$(\pi_b,\pi_w,\ell,\ell^*)$, with $\pi_b=\Pi(\phi,1) = \Pi_+(\phi)$,
$\pi_w=\Pi(\sigma \phi,2) = \Pi_+(H_- \sigma \phi)$, and
$\ell^*(\phi)$ if applicable,
\be
\Psi_\Lambda(\pi_b,\pi_w,\ell,\ell^*)
=
\Psi_\Lambda(\pi_w,\pi_b,\ell,\ell^*)
=
\Psi_\Lambda(R^{-1} \pi_b,R \pi_w,\ell,\ell^*)
\ef.
\ee
\end{theorem}

\section{The Razumov--Stroganov correspondence for the Scattering
  Matrix}
\label{main-section} 

\noindent
This section is mainly devoted to the proof of Theorem \ref{main},
that relates the solution of the scattering equation, on the $O(1)$
Dense Loop Model side, to the enumeration of FPL's performed according
to the refinement position, and the link pattern evaluated with the
map~$\Pi_b$.

\subsection{The correspondence for the enumerations according to $\Pi_b$}

\begin{figure}[!bt]
\[
\setlength{\unitlength}{3.75pt}
\begin{picture}(87,57.5)(-7,0)
\put(.5,22.8){\includegraphics[scale=.5]{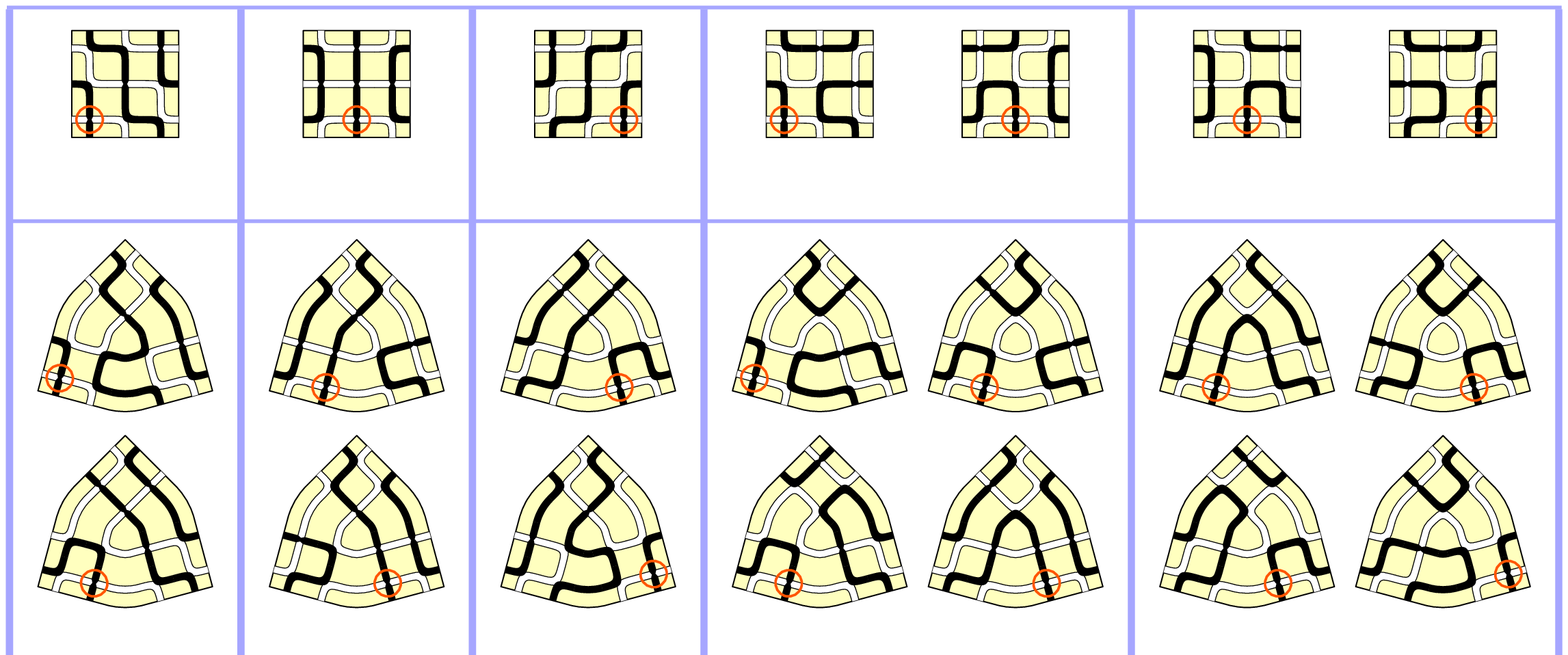}}
\put(-1.5,46.3){\makebox[0pt][r]{$\ket{\Psi}$:}}
\put(-1.5,20){\makebox[0pt][r]{$(1+t) \ket{\Psi}$:}}
\put(1,0){\makebox[0pt][r]{$(1-t)\beo_1 \ket{\Psi}$:}}
\put(6,0){\makebox[0pt][c]{$0$}}
\put(19,0){\makebox[0pt][c]{$0$}}
\put(32,0){\makebox[0pt][c]{$(1-t)(1+t+t^2)$}}
\put(50.5,0){\makebox[0pt][c]{$0$}}
\put(74.5,0){\makebox[0pt][c]{$(1-t)(1+t)^2$}}
\put(1,4){\makebox[0pt][r]{$-R\, \ket{\Psi}$:}}
\put(6,4){\makebox[0pt][c]{$-t$}}
\put(19,4){\makebox[0pt][c]{$-t^2$}}
\put(32,4){\makebox[0pt][c]{$-1$}}
\put(50.5,4){\makebox[0pt][c]{$-t(1+t)$}}
\put(74.5,4){\makebox[0pt][c]{$-(1+t)$}}
\put(1,8){\makebox[0pt][r]{$t \ket{\Psi}$:}}
\put(6,8){\makebox[0pt][c]{$t$}}
\put(19,8){\makebox[0pt][c]{$t^2$}}
\put(32,8){\makebox[0pt][c]{$t^3$}}
\put(50.5,8){\makebox[0pt][c]{$t(1+t)$}}
\put(74.5,8){\makebox[0pt][c]{$t^2(1+t)$}}
\put(6,20){\makebox[0pt][c]{$1+t$}}
\put(19,20){\makebox[0pt][c]{$t(1+t)$}}
\put(32,20){\makebox[0pt][c]{$t^2(1+t)$}}
\put(50.5,20){\makebox[0pt][c]{$(1+t)^2$}}
\put(74.5,20){\makebox[0pt][c]{$t(1+t)^2$}}
\put(6,46.3){\makebox[0pt][c]{$1$}}
\put(19,46.3){\makebox[0pt][c]{$t$}}
\put(32,46.3){\makebox[0pt][c]{$t^2$}}
\put(50.5,46.3){\makebox[0pt][c]{$(1+t)$}}
\put(74.5,46.3){\makebox[0pt][c]{$t(1+t)$}}
\end{picture}
\]
\caption{\label{fig.thetheo}Top: the 7 FPL on the $3\times 3$ square,
  collected according to the enumeration $\Pi_b$. Middle: analogous
  enumeration of the 14 FPL on the domain $\Lambda=(4,4;0,0,0,2)$.
  The enumerations are proportional, by a factor $1+t$. Bottom rows:
  the three summands involved in the scattering equation. These rows
  sum up to zero, so the scattering equation is satisfied by the
  enumeration vectors.}
\end{figure}

\noindent
In the previous sections we gathered all the ingredients required for
stating and proving our main result.  Assume that a dihedral domain
$\Lambda$ is given, and fix a reference side $r$ on the boundary, of
length $L$.
Call $h(\phi) \in \{1,\ldots,L\}$ the function associating to a FPL
configuration $\phi$ its refinement position along $r$, counted in
counter-clockwise order (starting from the reference corner).
Define the vector $\ket{\Psi_{\Lambda}(t)}$
whose components $\Psi_{\Lambda}(t;\pi)$ are enumerations of FPL in
$\fpl_{b}(\Lambda)$, with link pattern associated through the function
$\Pi_b$, and weighted according to the refinement position, i.e.
\be
\ket{\Psi_{\Lambda}(t)} := \sum_{\pi \in \LP(\Lambda)} \Psi_{\Lambda}(t;\pi)
\ket{\pi}
=
\sum_{\phi \in \fpl_{b}(\Lambda)}
t^{h(\phi)-1} 
\ket{\Pi_b(\phi)}
\ef.
\ee 
\begin{minipage}{\textwidth}
This section is devoted to the proof of the following
\begin{theorem}
\label{main}
The vector $\ket{\Psi_{\Lambda}(t)}$ defined above
satisfies the scattering equation at $i=1$, equation
(\ref{eq.qKZfund}), i.e.
\be
X_1(t)\ket{\Psi_{\Lambda}(t)}
= R\, \ket{\Psi_{\Lambda}(t)}
\ef.
\ee
In other words, for all choices of dihedral domain $\Lambda$, and
reference side $r$ of $\Lambda$,
\be
\label{eq.KinTh}
\ket{\Psi_{\Lambda}(t)}
=
K_{\Lambda}(t)\;
\ket{\Psi^{(1)}_{O(1)}(t)}
\ef,
\ee
where $K_{\Lambda}(t)$ is a polynomial with positive integer
coefficients, depending on $\Lambda$ and $r$, and
$\ket{\Psi^{(1)}_{O(1)}(t)}$ is the solution for the representation on
the space $\LP(\Lambda)$.
\end{theorem}
\noindent
An illustration of this theorem is given in Figure~\ref{fig.thetheo}.
\end{minipage}

For $\ast=+,-,b,w$, call $\fpl_{\ast}^{[i]}(\Lambda)$ the restriction
of $\fpl_{\ast}(\Lambda)$ to configurations such that $h(\phi)=i$. Call
$\ket{\Psi_{\Lambda}^{[i]}}$ the vector whose 
components are enumerations of $\fpl_{b}^{[i]}(\Lambda)$, using
$\Pi_b$. We can write $\ket{\Psi_{\Lambda}(t)}$ as
\be
\label{def-psi(t)}
\ket{\Psi_{\Lambda}(t)}
= 
\sum_{i=1}^{L}
t^{i-1}
\ket{\Psi_{\Lambda}^{[i]}}
=
\sum_{\pi \in \LP(\Lambda)}
t^{i-1}
\Psi_{\Lambda}^{[i]}(\pi) 
\ket{\pi}
=
\sum_{i=1}^{L}
t^{i-1}
\!\!\!
\sum_{\phi \in \fpl_{b}^{[i]}(\Lambda)}
\!\!\!
\ket{\Pi_b(\phi)}
\ef.
\ee
As a consequence of Proposition \ref{prop.equiqKZ}, the proof of
Theorem \ref{main} splits into the two lemmas
\begin{lemma}
\label{main-lemma1}
The vectors $\ket{\Psi_{\Lambda}^{[i]}}$ satisfy
\be
\label{main-theo-eq}
(\beo_1 - R\, \beo_N)
\ket{\Psi_{\Lambda}^{[i]}}=0
\ef;
\ee
\end{lemma}

\noindent
\begin{lemma}
\label{main-lemma2}
For any $\pi\in \LP(\Lambda)$ such that $1 \nsimX 2$ we have
\be
t \Psi_{\Lambda}(t;\pi) = \Psi_{\Lambda}(t,R^{-1} \pi)
\ef.
\ee
\end{lemma}
\noindent
We prove these lemmas in the remainder of this subsection.
Before this, we need some notations.  In Section
\ref{Temperley-section} we introduced a space
$\mathbb{C}^{\LP(\Lambda)}$, in order to conveniently encode the
action of Temperley--Lieb Algebra through linear operators.  Later
on, we have introduced maps $\Pi(\phi,v)$ from the set $\fpl(\Lambda)$
(or variants) to $\LP(\Lambda)$.  Here we will also consider the
spaces $\mathbb{C}^{\fpl_{\ast}(\Lambda)}$, with $\ast=+,-,b,w$, for
which a canonical basis is given by the configurations, $\kket{\phi}$,
in order to promote also the maps $\Pi(\phi,v)$ to linear operators.
Note that, for easiness of notation, we use a different graphical
representation for vectors in $\mathbb{C}^{\LP(\Lambda)}$, such as
$\ket{\pi}$, and vectors in $\mathbb{C}^{\fpl(\Lambda)}$, such as
$\kket{\phi}$.

With abuse of notation, we will still call $\Pi_+$ and $\Pi_b$ the
linear operators associated to the maps
$\Pi_+:\fpl_+(\Lambda)\rightarrow\LP(\Lambda)$ and
$\Pi_b:\fpl_b(\Lambda)\rightarrow \LP(\Lambda)$,
whose action on basis elements is $\Pi_{+}\kket{\phi} =
\ket{\Pi_{+}(\phi)}$ and $\Pi_{b}\kket{\phi} = \ket{\Pi_{b}(\phi)}$.
Define the vector associated to the collection of
configurations in $\fpl_b(\Lambda)$, weighted with their refinement
position, as
\be
\kket{s_{\Lambda}(t)}
:= 
\sum_{\phi \in \fpl_{b}(\Lambda)} 
t^{h(\phi)-1}
\kket{\phi}
\ef.
\ee
The vector $\ket{\Psi_{\Lambda}(t)}$ is thus the image of 
$\kket{s_{\Lambda}(t)}$ under the map $\Pi_{b}$,
$\ket{\Psi_{\Lambda}(t)}= \Pi_{b} \kket{s_{\Lambda}(t)}$.
Similarly, calling
\be
\kket{s_\Lambda^{[i]}}
:=
\sum_{\phi\in \fpl_b^{[i]}(\Lambda) }
\kket{\phi}
\ef,
\ee
we have
$\ket{\Psi_{\Lambda}^{[i]}}
= \Pi_{b} \kket{s_{\Lambda}^{[i]}}$.

\bigskip

\proofof{Lemma \ref{main-lemma1}} This lemma was already proven in
\cite{usRS} (Proposition 4.4, equation (61)), when the domain
$\Lambda$ was a $n \times n$ square. Actually, the proof easily extends to
the general class of domains $\Lambda$ that we described in Section
\ref{ssec.multic}, which is also the generality of the treatment of
Wieland gyration discussed in Section 3 of~\cite{usRS}.

We want to give here a slightly reformulated proof. First, we define
two involutions $\tilde \beo_1$ and $\tilde \beo_N$ on
$\fpl(\Lambda)$.  These operators are conjugated versions of the local
gyration $H_{\gamma}$, for $\gamma_1$ and $\gamma_N$ certain cycles in
the graph, adjacent to the refinement position.  If $\phi \in \fpl_+$, 
$\gamma_1$ and $\gamma_N$ are the faces at the top-right and top-left
corner of the refinement $h(\phi)$, respectively.  If $\phi \in
\fpl_-$, right and left are interchanged.

The operator $\tilde \beo_1$ flips the cycle $\gamma_1$, if it has
length 4 and it is coloured in an alternating way, i.e.~of the form
$\sqqV$ or $\sqqH$\;\footnote{Or if it has length 2, is punctured, and
  is coloured in an alternating way, but one can see that non-trivial
  domains never have such a face with adjacent to the boundary.},
otherwise it leaves the FPL unchanged. The operator $\tilde \beo_N$
has an analogous action on $\gamma_N$. The local configurations which
are not kept fixed are
\[
\setlength{\unitlength}{81pt}
\begin{picture}(4.45,.44)(0,0)
\put(0.066,0.1){$\tilde \beo_1:$}
\put(2.366,0.1){$\tilde \beo_N:$}
\put(0.3,0){\includegraphics[scale=0.9]{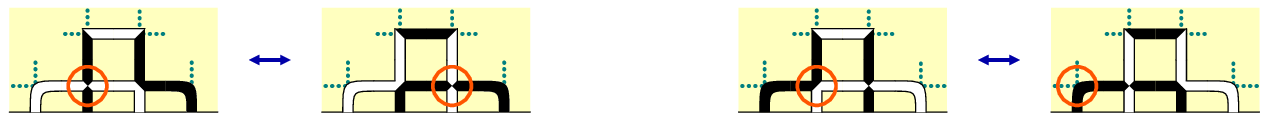}}
\end{picture}
\]
If $\phi \in \fpl_b(\Lambda)$, the operators $\tilde \beo_1$ and
$\tilde \beo_N$ act analogously to the Temperley--Lieb generators
$\beo_1$ and $\beo_N$, in the sense that it holds
\begin{subequations}
\label{eqs.comm-tildee-pi}
\begin{align}
\label{comm-tildee-pi2}
\Pi(\tilde \beo_N \phi, h(\phi))
&= 
\beo_N
\Pi(\phi,h(\phi))= 
\beo_N \Pi_b(\phi)
\ef;
\\
\label{comm-tildee-pi}
\Pi(\tilde \beo_1 \phi, h(\phi))
&= 
\beo_1 
\Pi(\phi,h(\phi))= 
\beo_1 
\Pi_b(\phi)
\ef.
\end{align}
\end{subequations}
Consider the state $\kket{s_\Lambda^{[i]}}$ (so $i=h(\phi)$ for
all contributing configurations). We claim that
\be
\label{gyr-rels}
H\,
\tilde \beo_N
\kket{s_\Lambda^{[i]}}
=
\sigma
\kket{s_\Lambda^{[i]}}
=
H^{-1}
\tilde \beo_1
\kket{s_\Lambda^{[i]}}
\ef.
\ee
In other words, the three maps $\sigma$, $H\, \tilde \beo_N$ and
$H^{-1} \tilde \beo_1$ are bijections from $\fpl_b^{[i]}(\Lambda)$ to
$\fpl_w^{[i]}(\Lambda)$ (this is clear for $\sigma$). It is obvious,
from the fact that $\sigma$, $H$, $\tilde \beo_1$ and $\tilde \beo_N$
are all invertible, that they are bijections on $\fpl(\Lambda)$, and
the only fact that needs to be checked is that, in both non-trivial
cases, the image of a FPL with refinement position $i$, and black, is
a FPL with refinement position $i$, and white.  Figure \ref{fig.lem41}
shows this fact for $H\, \tilde \beo_N$.

\begin{figure}[!bt]
\[
\setlength{\unitlength}{10pt}
\begin{picture}(31,11)
\put(0,0){\includegraphics[scale=2]{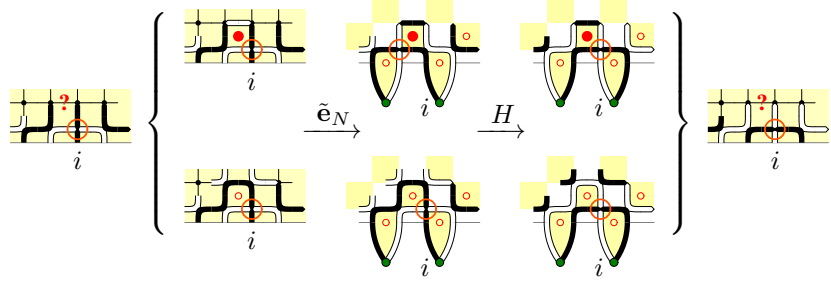}}
\put(2.85,4){$i$}
\put(28.95,4){$i$}
\put(9.35,7){$i$}
\put(9.35,1){$i$}
\put(15.85,6){$i$}
\put(15.85,0){$i$}
\put(22.35,6){$i$}
\put(22.35,0){$i$}
\put(5.5,5.5){$\displaystyle{\left\{\rule{0pt}{45pt}\right.}$}
\put(11.4,5.2){$\xrightarrow{~\displaystyle{\tilde \beo_N}}$}
\put(17.9,5.2){$\xrightarrow{~\displaystyle{H}}$}
\put(25.,5.5){$\displaystyle{\left.\rule{0pt}{45pt}\right\}}$}
\end{picture}
\]
\caption{\label{fig.lem41}Analysis of the map $H \tilde \beo_N$ on
  $\fpl_b^{[i]}(\Lambda)$.
  The state $\kket{s_\Lambda^{[i]}}$ is decomposed according to the
  colour of the edge marked with `?'. The action of $\tilde \beo_N$ is
  non-trivial if and only if this edge is white. Then, we can analyse
  the action of $H$, by producing the graph $\Lambda_+$. The result is
  the decomposition, according to the colour of the edge marked with
  `?', of the state $\sigma \kket{s_\Lambda^{[i]}}$.}
\end{figure}

Equations (\ref{eqs.comm-tildee-pi}) allow to deduce
\begin{subequations}
\label{eqs.5435765x}
\begin{align}
\begin{split}
\beo_N \ket{\Psi_{\Lambda}^{[i]}}
&=
\beo_N \Pi_b \kket{s_{\Lambda}^{[i]}}
=
\Pi(\cdot,i)
\,
\tilde \beo_N(i)
\,
\kket{s_{\Lambda}^{[i]}}
\ef.
\end{split}
\label{eq.5435765}
\\
\label{eq.5435765b}
\begin{split}
\beo_1 \ket{\Psi_{\Lambda}^{[i]}}
&=
\beo_1 \Pi_b \kket{s_{\Lambda}^{[i]}}
=
\Pi(\cdot,i)
\,
\tilde \beo_1(i)
\,
\kket{s_{\Lambda}^{[i]}}
\ef.
\end{split}
\end{align}
\end{subequations}
We can now compare the two left hand sides of (\ref{eqs.5435765x}),
using (\ref{gyr-rels}) and (\ref{wiel-gyr-fpl})
\be
\begin{split}
\beo_1 \ket{\Psi_{\Lambda}^{[i]}}
&=
\Pi(\cdot,i)
\,
H^2
\tilde \beo_N(i)
\,
\kket{s_{\Lambda}^{[i]}}
=
R
\,
\Pi(\cdot,i)
\,
\tilde \beo_N(i)
\,
\kket{s_{\Lambda}^{[i]}}
=
R\, \beo_N \ket{\Psi_{\Lambda}^{[i]}}
\ef;
\end{split}
\ee
as was to be proven.
\qed

\bigskip

\begin{figure}[!b]
\[
\setlength{\unitlength}{10pt}
\begin{picture}(25,6.5)
\put(0,0){\includegraphics[scale=2]{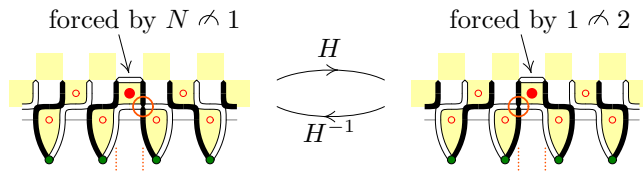}}
\put(12.5,4.4){\makebox[0pt][c]{$H$}}
\put(12.5,1.2){\makebox[0pt][c]{$H^{-1}$}}
\put(2,5.5){forced by $N \nsimX 1$}
\put(17,5.5){forced by $1 \nsimX 2$}
\end{picture}
\]
\caption{\label{fig.statesXY}Analysis of Wieland half-gyration in a
  neighbourhood of the refinement position, required in the proof of
  Lemma \ref{main-lemma2}.  }
\end{figure}

\proofof{Lemma \ref{main-lemma2}} Consider a pattern $\pi\in
\LP(\Lambda)$ such that $1 \nsimX 2$. We claim that, for every $\phi
\in \fpl_{b}(\Lambda)$ such that $\Pi_b(\phi)=\pi$ and $h(\phi)=h$,
the configuration $\phi'=H^{-1}(\phi)$ has $\Pi_b(\phi')=R^{-1} \pi$
and $h(\phi')=h+1$. In particular it has $N \nsimX 1$. This is easily
seen from the local properties of gyration, in a neighbourhood of the
refinement position (see Figure \ref{fig.statesXY}, right to
left). Analogously, for a pattern $\pi\in \LP(\Lambda)$ such that $N
\nsimX 1$, for every $\phi \in \fpl_{b}(\Lambda)$ such that
$\Pi_b(\phi)=\pi$ and $h(\phi)=h$, the configuration $\phi'=H(\phi)$
has $\Pi_b(\phi')=R\, \pi$ and $h(\phi')=h-1$ (see Figure
\ref{fig.statesXY}, left to right).

Thus, $H$ induces a bijection between the sets of configurations in
$\fpl_{b}(\Lambda)$ with link pattern $\pi$ and $R^{-1}\pi$, given
that $1 \nsimX 2$ in $\pi$ and thus $N \nsimX 1$ in $R^{-1} \pi$. As
this bijection acts in a simple way on the refinement position, it
induces a relation at the level of the weighted generating functions,
$\Psi_{\Lambda}(t;\pi)$ and $\Psi_{\Lambda}(t;R^{-1}\pi)$, which is
exactly the statement of the lemma.  \qed

\subsection{Special families of dihedral domains}
\label{numerics}

\noindent
In the definition of the vector $\ket{\Psi_\Lambda(t)}$ there is an
implicit dependence from the reference side.  A generic dihedral
domain has several (at most $4$) possible choices of reference side $r$,
which, except for the case of the $n \times n$ square and a few other
symmetric cases, are in general not equivalent.

In this section we analyse this aspect, so, only within this section,
we adopt a notation that makes explicit the dependence on the
reference side $r$, by writing $\ket{\Psi_{{\Lambda,r}}(t)}$ for the
vector, and $K_{{\Lambda,r}}(t)$ for the proportionality factor with
the minimal solution of the scattering equation, i.e.
\be
\ket{\Psi_{{\Lambda,r}}(t)} =
K_{{\Lambda,r}}(t)
\,
\ket{\Psi_{O(1)}^{(1)}(t)}
\ef.
\ee
Interestingly, for all the three possible realisations on the
correspondence, $\LP(\Lambda) = \LP(2n)$, $\LP^*(2n)$ or
$\LP^*(2n-1)$, there exists a dihedral domain (and a reference side)
realising the minimal degree, i.e.\ such that $K_{{\Lambda,r}}(t)$
is a constant, and in fact can be set to $1$.

For $\LP(\Lambda) = \LP(2n)$, the realisation with minimal degree
corresponds to the $n \times n$ square domain, while for $\LP(\Lambda)
= \LP^*(N)$ (either even or odd) the realisation with minimal degree
corresponds to Half-turn symmetric ASM of size $N$ (refer back to
Figure~\ref{fig.symclass}).  Indeed, ordinary ASM's on the square of
side $n$ have refinement positions in the range $\{1,\ldots,n\}$, thus
the corresponding polynomials $\Psi_{\Lambda}(t;\pi)$ have degree at
most $n-1$. Similarly, HTASM's of side $N$ have refinement positions
in the range $\{1,\ldots,N\}$, thus the corresponding polynomials
$\Psi_{\Lambda}(t;\pi)$ have degree at most $N-1$.  As observed at the
end of Section \ref{ssec.multic}, both ASM and HTASM are special cases
of dihedral domains, and the latter allows for the presence of a
puncture.  Furthermore, it is easily seen that for each rainbow
pattern in $\LP(2n)$ there exists a unique ASM of size $n$, and that
for each rainbow pattern in $\LP^*(N)$ there exists a unique HTASM of
size $N$.

As we have proven that ASM's on the square present refined dihedral
Razumov--Stroganov correspondence on the link-pattern space $\LP(2n)$,
and that HTASM's present refined dihedral Razumov--Stroganov
correspondence on the link-pattern space $\LP^*(N)$, the reasonings of
the previous paragraph imply upper bounds on the degree of
$\Psi_{O(1)}^{(1)}(t)$, in the two cases $\LP(2n)$ and $\LP^*(N)$. As
these bounds match the lower bounds in Corollary \ref{coroll.degLB},
we get the following

\begin{proposition}
\label{prop.degT}
For $\LP(2n)$, there exist polynomial solutions
$\ket{\Psi^{(i)}(t)}$ 
of the scattering equation
(\ref{eq.qKZfund}), such that the components
$\Psi^{(i)}(t;\pi)$ 
are polynomials in $t$ with positive
integer coefficients, the maximal (over $\pi$'s) degree is $n-1$, and
$\Psi^{(i)}(t;\pi) = t^j$ for some $j$ if $\pi$ is a rainbow
pattern. For $\LP^*(N)$, the same facts do hold, with maximal degree
$N-1$.
\end{proposition}

\bigskip
\noindent
For a generic dihedral domain $\Lambda$, the space $\LP(\Lambda)$ is
independent from the reference side $r$, and thus also the DLM vector
$\ket{\Psi_{O(1)}^{(1)}(t)}$ is the same. However, since the external
lines are possibly of different lengths, and the range of allowed
possible refinement positions, i.e. indices $h$ such that
$\ket{\Psi_{{\Lambda,r}}^{[h]}}\neq 0$, may be different for different
sides $r$, the vectors $\ket{\Psi_{{\Lambda,r}}(t)}$, for a given
$\Lambda$ and the various possible choices of $r$, are in general
distinct polynomials, and even of different degree. Our result implies
that this difference arises uniquely at the level of the overall
factor $K_{{\Lambda,r}}(t)$.



A sub-family of dihedral domains in which this feature can be seen
explicitly is the one with three corners (this is the maximum allowed
number, besides the `classical' $n \times n$ square).  We can set
$a_1=a_2=a_3=0$, and $a_4>0$, without loss of generality.  These
domains also minimize the number of `crossing bundles' (again, besides
the case of the square): we have three bundles crossing each other, of
width\footnote{A bundle has width $k$ if it is constituted of $k$
  parallel lines.}  $a_4$, $L_x-a_4$ and $L_y-a_4$.

\begin{figure}
\[
\setlength{\unitlength}{10pt}
\begin{picture}(19,16.5)
\put(0.7,0.7){\includegraphics[scale=.666]{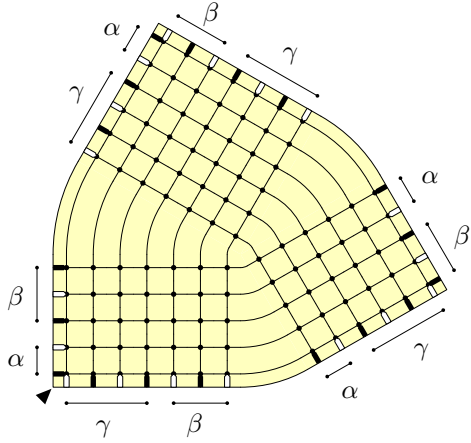}}
\put(3.3,0){$\gamma$}
\put(6.7,0){$\beta$}
\put(12.5,1.3){$\alpha$}
\put(15.2,2.9){$\gamma$}
\put(16.7,7.1){$\beta$}
\put(15.5,9.4){$\alpha$}
\put(10.3,14.1){$\gamma$}
\put(7.3,15.5){$\beta$}
\put(3.6,15){$\alpha$}
\put(2.3,12.6){$\gamma$}
\put(0,4.8){$\beta$}
\put(0,2.5){$\alpha$}
\end{picture}
\]
\caption{\label{fig.trefasci}A three-bundle domain
  $T(\alpha,\beta,\gamma)$, with $(\alpha,\beta,\gamma)=(2,3,4)$. The
  reference side has length $\alpha+2\gamma+\beta=13$.}
\end{figure}

\begin{figure}
\[
\setlength{\unitlength}{14.2pt}
\begin{picture}(12.6,9.2)
\put(0.45,0.66){\includegraphics[scale=1.38]{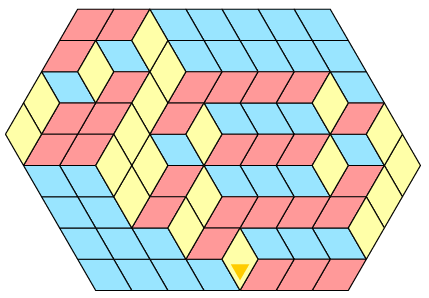}}
\put(0.8,7.3){$\alpha$}
\put(1.1,2.2){$\beta$}
\put(6,9.4){$\gamma$}
\put(3,0){$1$}
\put(4,0){$2$}
\put(5,0){$3$}
\put(6,0){$4$}
\put(7,0){${\bf 5}$}
\put(8,0){$6$}
\put(9,0){$7$}
\put(10,0){$8=\gamma+1$}
\end{picture}
\quad
\begin{picture}(12.6,6)
\put(0.45,0.66){\includegraphics[scale=1.38]{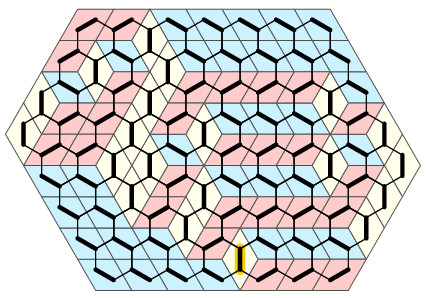}}
\put(0.8,7.3){$\alpha$}
\put(1.1,2.2){$\beta$}
\put(6,9.4){$\gamma$}
\put(3,0){$1$}
\put(4,0){$2$}
\put(5,0){$3$}
\put(6,0){$4$}
\put(7,0){${\bf 5}$}
\put(8,0){$6$}
\put(9,0){$7$}
\put(10,0){$8=\gamma+1$}
\end{picture}
\]
\caption{\label{fig.exMacMah}Left: example of weighted rhombus tiling
  corresponding to the expression in (\ref{fact-tr}). Right: the
  associated dimer covering of a region of the honeycomb lattice.}
\end{figure}

A useful and more symmetric parametrisation consists in setting
$L_x=\alpha+2\beta + \gamma$, $L_y=\alpha+\beta + 2\gamma$ and
$a_4=\beta+\gamma$, so that the bundles have widths $\alpha+\beta$,
$\alpha+\gamma$ and $\beta+\gamma$.  We call such a domain a
\emph{three-bundle domain}, or \emph{triangoloid},
$T(\alpha,\beta,\gamma)$ (see Figure \ref{fig.trefasci}).  There are
three possible choices of reference sides, $r(\alpha)$, $r(\beta)$ and
$r(\gamma)$, of lengths
$2\alpha+\beta+\gamma$, $\alpha+2\beta+\gamma$ and
$\alpha+\beta+2\gamma$, respectively.
The number of external edges of each colour is
$2n=2(\alpha+\beta+\gamma)$ (we use $n$ as a synonim of
$\alpha+\beta+\gamma$ in what follows).  Since
$T(\alpha,\beta,\gamma)$ has no faces with less than three edges, the
correspondence holds with the non-punctured link patterns,
$\LP(\Lambda)=\LP(2n)$.
%
We thus have
$
\ket{\Psi_{T(\alpha,\beta,\gamma),r(\gamma)}(t)} =
K_{T(\alpha,\beta,\gamma),r(\gamma)}(t)
\; \ket{\Psi_{\LP(2n)}(t)} 
$,
or, with a more compact \emph{ad hoc} notation,
\be
\ket{\Psi_{\alpha,\beta|\gamma}(t)} =
K_{\alpha,\beta|\gamma}(t)
\;
\ket{\Psi_{\LP(2n)}(t)} 
\ef.
\ee
Simple reasonings (of reflection symmetry) show that
$K_{\alpha,\beta|\gamma}(t)$ is symmetric under interchange of
$\alpha$ and $\beta$, and $t \to 1/t$ (up to an appropriate rescaling
by $t^{\gamma}$, because $\gamma=L-n$). Here we report a
formula for $K_{\alpha,\beta|\gamma}(t)$, without providing a
proof\;\footnote{The proof of equation (\ref{fact-tr}) goes through
  computing the partition function of the $6$-Vertex Model on the
  triangoloid $T(\alpha,\beta,\gamma)$, with \emph{generic} spectral
  parameters $z_i$ and at $q=e^{2\pi i/3}$.}
%
\be
\label{fact-tr}
\begin{split}
K_{\alpha,\beta|\gamma}(t)
&=
\det_{1\leq i,j\leq \gamma}
\left(\binom{\alpha+\beta-1}{\beta-i+j}
+t\,\binom{\alpha+\beta-1}{\beta-i+j-1} \right)  
\\
&
=
\frac{\Delta_{\alpha} \Delta_{\beta}
\Delta_{\gamma+1}
\Delta_{\alpha+\beta+\gamma-1}}
{\Delta_{\alpha+\beta-1}
\Delta_{\alpha+\gamma}\Delta_{\beta+\gamma}
}
\;
\sum_{i=0}^{\gamma} t^i  
\binom{\beta-1+i}{i}
\binom{\alpha-1+\gamma-i}{\gamma-i}
\ef,
\end{split}
\ee
where $\Delta_n = \prod_{k=0}^{n-1} k!$.
The function $K_{\alpha,\beta|\gamma}(t)$ has also a combinatorial
interpretation in terms of the weighted enumerations of rhombus
tilings of a hexagon with sides $(\alpha,\beta,\gamma)$. As well as
ASM's, rhombus tilings of a hexagon have a natural notion of
refinement position w.r.t.\ a given side (defined as the position of
the only lozenge, adjacent to the side, with the longer diagonal
orthogonal to the side). In our case, a tiling has a weight $t^{h-1}$,
for $h$ the refinement position w.r.t.\ one of the two sides of length
$\gamma$ (see Figure~\ref{fig.exMacMah}).

This connection between lozenge tilings of a hexagon and the factor
$K_{\Lambda,r}$ for the Razumov--Stroganov correspondence on
three-bundle domains has also a direct combinatorial interpretation
through an analysis of frozen regions in configurations with rainbow
link patterns.

\subsection{Specialisations to $t=1$ and $t=0$}

\noindent
As explained in Section \ref{exch-section}, there are two relavant
specialisations of the parameter $t$. When $t=1$, the vector
$\ket{\Psi_{\Lambda}(1)}$ is the ground state of the Hamiltonian of
the $O(1)$ Dense Loop Model, equation (\ref{hamil-O(1)}). This means
that we have a realisation of the Razumov--Stroganov correspondence
involving a \emph{different} class of FPL's on a domain $\Lambda$,
namely $\fpl_b(\Lambda)$ and a \emph{different} function associating a
link pattern, namely $\Pi_b(\phi)$.  One way to obtain from this
result the usual Razumov--Stroganov correspondence, involving
$\fpl_+(\Lambda)$ and $\Pi_+(\phi)$, will be presented in
Section~\ref{bij-section}.

Another way to recover the usual correspondence is by looking at the
specialisations $t \to 0$ (or at the symmetric limit $t \to \infty$).
In Section \ref{exch-section} we have recalled a result of Di
Francesco and Zinn-Justin \cite{pdf-pzj1}, which implies that the
specialisation $t=0$ of a solution of the scattering equation
(\ref{eq.qKZfund}), at $i=1$, is equal to the image of the ground
state of the Hamiltonian (\ref{hamil-O(1)}) of size $N-2$, under the
map that inserts a short arc $(N,1)$ (see equation
(\ref{recursion-paul-philippe})).
A correspondence is found if we can show that the analogous feature 
arises also on the FPL side.

Consider a domain $\Lambda=(L_x,L_y;0,a_2,a_3,a_4)$, such that all the
parameters $a_2$, $a_3$ and $a_4$ are at least 2 (in particular,
$\Lambda$ has a single corner).  In such a domain it is easily seen
that $\ket{\Psi_{\Lambda}^{[1]}} = 0$, while, as we now deduce, the
expression $\ket{\Psi_{\Lambda}^{[2]}}=\lim_{t \to 0} \frac{1}{t}
\ket{\Psi_{\Lambda}(t)}$ can be related to the ordinary enumeration of
FPL on a suitable reduced domain $\Lambda'$.

\begin{figure}
\[
\setlength{\unitlength}{12.5pt}
\begin{picture}(24.5,10)
\put(0,0){\includegraphics[scale=1.25]{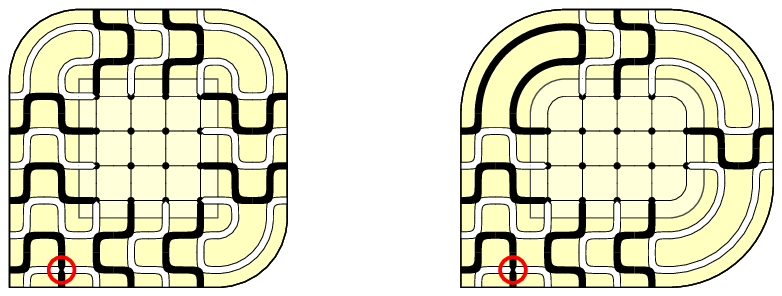}}
\put(3,0.5){$1$}
\put(5,0.5){$2$}
\put(7,0.5){$3\ \ldots$}
\put(1.5,1.5){\makebox[0pt][r]{$N$}}
\put(1.5,3.5){\makebox[0pt][r]{$N-1$}}
\put(16,0.5){$1$}
\put(18,0.5){$2$}
\put(20,0.5){$3\ \ldots$}
\put(14.5,1.5){\makebox[0pt][r]{$N$}}
\put(14.5,3.5){\makebox[0pt][r]{$N-1$}}
\end{picture}
\]
\caption{\label{fig.tto0}Frozen region induced by imposing that the
  refinement position is 2, on a domain
  $\Lambda=(L_x,L_y;0,a_2,a_3,a_4)$ such that all the parameters
  $a_2$, $a_3$ and $a_4$ are at least 2. The light-yellow unfrozen
  region is the associated reduced domain $\Lambda'$. Left: $\Lambda'$
  corresponds to ordinary ASM's, and has the maximal number of
  corners, 4. Right: $\Lambda'$ has the minimal number of corners, 1.
}
\end{figure}

More precisely, within $\ket{\Psi_{\Lambda}^{[2]}}$ the only
non-vanishing components have $1\simX N$, as this connection is
implied by the tile at the refinement position and the induced frozen
region (and the puncture, if any, is within the unfrozen region and
thus not surrounded by this short arc). The unfrozen region
corresponds to the domain
$\Lambda'=(L_x-4,L_y-4;0,a_2-2,a_3-2,a_4-2)$, and has induced
alternating boundary conditions, with a white external edge adjacent
to the reference corner. Thus the FPL's $\phi$ in
$\fpl_b^{[2]}(\Lambda)$ are in bijection with the FPL's $\phi'$ in
$\fpl_-(\Lambda')$. The external edges in $\Lambda'$ are connected to
the external edges in $\Lambda$, numbered from $2$ to $N-1$, in
counter-clockwise order starting from the reference corner, thus,
for $\phi$ and $\phi'$ configurations in bijection, the link patterns
$\Pi_-(\phi')$ and $\Pi_b(\phi)$ are easily related.
See also Figure \ref{fig.tto0}. This gives 
\be
\label{recursion-paul-philippe.FPL.2}
\Psi_{\Lambda}^{[2]}(
\setlength{\unitlength}{6pt}
\begin{picture}(8,2.8)
\put(0,0){\includegraphics[scale=1.2]{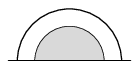}}
\put(0.6,-1){${}_{1}\rule{30pt}{0pt}{}_N$}
\put(3.5,0.6){$\pi$}
\end{picture}
)
=
\Psi_{\Lambda'}(1,\pi)
\ef,
\ee
and therefore, in view of equation (\ref{recursion-paul-philippe}),
$\ket{\Psi_{\Lambda'}(1)}$ is the ground state of the Hamiltonian
corresponding to the domain $\Lambda'$, which means that the
Razumov--Stroganov correspondence holds in the domain~$\Lambda'$.

As for every dihedral domain $\Lambda'=(L'_x,L'_y;0,a'_2,a'_3,a'_4)$
we can construct a dihedral domain
$\Lambda=(L'_x+4,L'_y+4;0,a'_2+2,a'_3+2,a'_4+2)$ with the required
properties, the reasoning above produces an alternative proof of the
original Razumov--Stroganov correspondence (on all dihedral domains,
and also in the punctured version).

\section{From the enumerations $\Psi_{\Lambda}(t;\pi)$ to the ordinary
  Razumov--Stroganov correspondence and Di Francesco's
  conjecture}\label{from-main-to-DF}

\noindent
At the end of the previous section we described a way of deriving the
ordinary Razumov--Stroganov correspondence from Theorem \ref{main},
that makes use of a deep and general result in \cite{pdf-pzj1} for
relating the limit $t \to 0$ of the solution of the scattering
equation to the ground state of the Hamiltonian at smaller size.

Here, in Section \ref{ssec.ph04proof} we provide a derivation of the
more general Di Francesco's conjecture in \cite{PdF04}, quickly
described in the introduction, and whose precise statement is reported
in the following Theorem \ref{df-conj}. At the light of the reasonings
in Section \ref{exch-section}, this also provides an alternative,
self-contained proof of the ordinary Razumov--Stroganov
correspondence.  Then, in Section \ref{bij-section} we give a third,
bijective derivation of the ordinary Razumov--Stroganov
correspondence.  Also this derivation does not rely on
\cite{pdf-pzj1}, and furthermore it doesn't make use of the results of
Section \ref{ssec.ph04proof}.

All these results follow from Theorem \ref{main}, and from an analysis
of the structure of the orbits
under the action of the half-gyration $H$, presented in Section
\ref{ssec.orbitstr}.

\subsection{Structure of the orbits under half-gyration}
\label{ssec.orbitstr}

\noindent
The \emph{orbit} associated to a configuration $\phi\in
\fpl(\Lambda)$, under the action of the half-gyration operator $H$, is
defined as the sequence $(\phi, \,H\, \phi, \,H^2\, \phi, \,H^3\,
\phi, \ldots, \,H^{p-1}\, \phi)$, where $p=p(\phi)$ is the smallest
positive integer such that $H^p \phi = \phi$ (i.e., the \emph{period}
of the orbit).

Remark that the configurations in the list are alternating in
$\fpl_+(\Lambda)$ and $\fpl_-(\Lambda)$, and in particular the period
must be even.  The precise periodicity of the various orbits is
immaterial at our purposes, and, in order to disentangle this element
from the analysis, we will mostly concentrate equivalently on the
periodically-repeated infinite sequences.

Thus, for some configuration $\phi \in \fpl(\Lambda)$, we define
$\cO(\phi)$ as the infinite sequence 
$\cO(\phi) := \{ \phi_t \} _{t \in \mathbb{Z}}$, determined
as $\phi_0 \equiv \phi$, and $\phi_{t} = H^t \phi_0$ for $t \neq 0$, 
i.e.
\[
\ldots \; \xrightarrow{H} \, 
\phi_{-2} \; \xrightarrow{H} \, 
\phi_{-1} \; \xrightarrow{H} \, 
\phi_0 
\; \xrightarrow{H} \, \phi_1
\; \xrightarrow{H} \, \phi_2
\; \xrightarrow{H} \, \ldots
\]
(recall that $H$ is invertible). Call $\cO^p$ a periodic portion of an
infinite orbit $\cO$, and $\cO^p_{\ast}$ the part of $\cO^p$ in
$\fpl_{\ast}(\Lambda)$, for $\ast=+,-,b,w$.

Recall that we defined above $h(\phi) \in \{1,\ldots,L\}$ as the
refinement position of $\phi$, i.e., position of the only $c$-type
tile along our reference side. Define also $d(\phi) \in \{a,b,c\}$ as
the direction taken at the refinement position, i.e., the kind of
tile, among $\{a,b,c\}$, for the tile `immediately above' the
refinement position (this is well-defined for all dihedral domains,
except for a few cases with size of order 1, like the $1 \times 1$
square). Note that the notion of $d(\phi)$ uses the fact that the
internal vertices have degree 4, but is independent from the number of
sides in the neighbouring faces, and is thus well defined also for
those dihedral domains in which some of the faces adjacent to the
reference side are triangles. As a shortcut, for a fixed orbit
$\cO(\phi)$, call $h_t = h(\phi_t)$ and $d_t = d(\phi_t)$.

We have the following useful lemma:
\begin{lemma}
\label{lem.orbitpattern}
If $\phi \in \fpl_+(\Lambda)$, the only possible local patterns
for the sequence $\{(h_t, d_t)\}_{t \in \mathbb{Z}}$ associated to its
orbit $\cO(\phi)$ are
\be
\label{h-patterns}
\begin{array}{r|rcccl|c|c|}
\cline{2-8}
&
\multicolumn{5}{|c|}{\textrm{local pattern}}
& \multicolumn{2}{|c|}{\textrm{condition}} \\
\cline{7-8}
\rule{0pt}{10.5pt}
&
\makebox[7pt][l]{$\cdots$} &h_{t-1} & h_t & h_{t+1} & \makebox[7pt][r]{$\cdots$}
                       & d_t & h_t-t \\
\cline{2-8}
\rule{0pt}{11pt}
&\makebox[7pt][l]{$\cdots$} &h   &  h &  h+1& \makebox[7pt][r]{$\cdots$}  & a & \textrm{even} \\
\phi \in \fpl_w(\Lambda)
\setlength{\unitlength}{8pt}
\begin{picture}(1,1)
\put(0.5,0){$\displaystyle{
\left\{ \rule{0pt}{16pt} \right.
}$}
\end{picture}
&\makebox[7pt][l]{$\cdots$} &h-1 &  h &  h+1& \makebox[7pt][r]{$\cdots$}  & c & \textrm{even} \\
&\makebox[7pt][l]{$\cdots$} &h-1 &  h &  h& \makebox[7pt][r]{$\cdots$}    & b & \textrm{even} \\
&\makebox[7pt][l]{$\cdots$} &h+1 &  h &  h& \makebox[7pt][r]{$\cdots$}    & a & \textrm{odd} \\
\phi \in \fpl_b(\Lambda)
\setlength{\unitlength}{8pt}
\begin{picture}(1,1)
\put(0.5,0){$\displaystyle{
\left\{ \rule{0pt}{16pt} \right.
}$}
\end{picture}
&\makebox[7pt][l]{$\cdots$} &h+1 &  h &  h-1& \makebox[7pt][r]{$\cdots$}  & c & \textrm{odd} \\
&\makebox[7pt][l]{$\cdots$} &h   &  h &  h-1& \makebox[7pt][r]{$\cdots$}  & b & \textrm{odd} \\
\cline{2-8}
\end{array}
\phantom{descending \rule{8pt}{0pt}}
\ee
If $\phi \in \fpl_-(\Lambda)$, even and odd are interchanged in the
column~$h_t-t$.
\end{lemma}
Another illustration of the table in (\ref{h-patterns}) is
given in Figure~\ref{fig.tabproofautom}.

\begin{figure}[!tb]
\[
\setlength{\unitlength}{25pt}
\begin{picture}(14.5,7.3)(0,0.5)
\put(0,0){\includegraphics[scale=2.5]{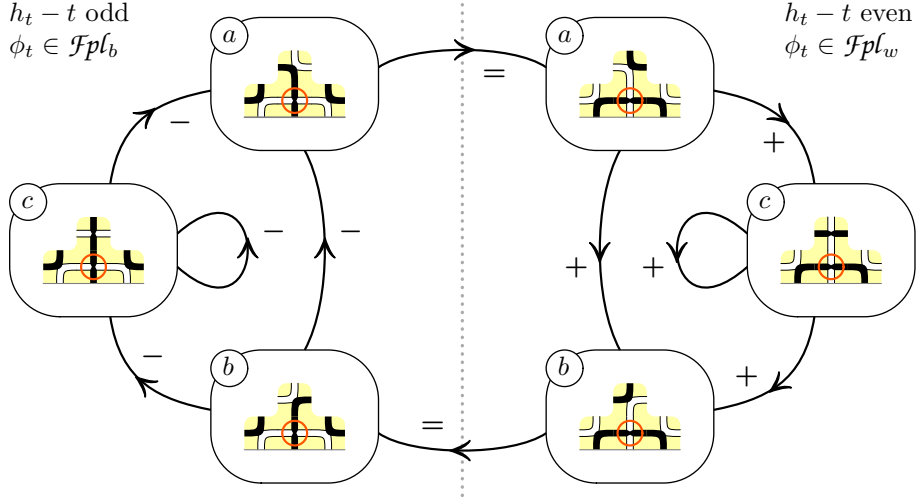}}
\put(0.5,7.2){\begin{minipage}{70pt}
$h_t-t$ odd\\
$\phi_t \in \fpl_b$
\end{minipage}}
\put(12.1,7.2){\begin{minipage}{70pt}
$h_t-t$ even\\
$\phi_t \in \fpl_w$
\end{minipage}}
\put(3.69,7.16){$a$}
\put(3.69,2.13){$b$}
\put(0.68,4.66){$c$}
\put(8.71,7.16){$a$}
\put(8.71,2.13){$b$}
\put(11.72,4.66){$c$}
\put(7.58,6.59){$\bm{=}$}
\put(6.66,1.27){$\bm{=}$}

\put(5.43,4.2){$\bm{-}$}
\put(4.28,4.2){$\bm{-}$}
\put(2.88,5.85){$\bm{-}$}
\put(2.47,2.31){$\bm{-}$}

\put(8.81,3.66){$\bm{+}$}
\put(9.96,3.66){$\bm{+}$}
\put(11.36,2.01){$\bm{+}$}
\put(11.77,5.55){$\bm{+}$}

\end{picture}
\]
\caption{\label{fig.tabproofautom}A different illustration of the
  statement of Lemma~\ref{lem.orbitpattern}. We show the possible
  transitions, in the set of six states $(d_t, h_t-t) \in \{a,b,c\} \times
  \{\textrm{even},\textrm{odd}\}$. The symbol in $\{-,=,+\}$ next to
  the arrow describes the value of $h_{t+1}-h_t$, in $\{-1,0,+1\}$
  respectively.}
\end{figure}

\proof The lemma is obtained by investigation of the action of $H$, at
the light of the frozen regions induced by the specialisation of
$d_t$. Remark that the parity of $h_t-t$ corresponds to the fact that
the half-gyration acting at time $t$ acts on the plaquette immediately
at the right or at the left of the refinement position.

\begin{figure}[!tb]
\[
\setlength{\unitlength}{25pt}
\begin{picture}(15,12)
\put(0,0){\includegraphics[scale=2.5]{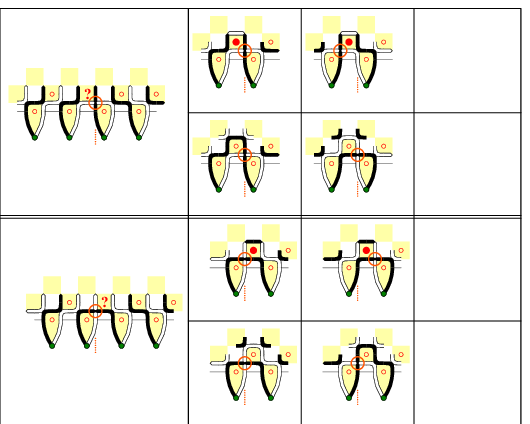}}
\put(0.5,7.2){\begin{minipage}{70pt}
$h_t-t$ odd\\
$\phi_t \in \fpl_b$
\end{minipage}}
\put(0.5,1.2){\begin{minipage}{70pt}
$h_t-t$ even\\
$\phi_t \in \fpl_w$
\end{minipage}}
\put(6,9.25){$d_t = b,c$}
\put(6,6.25){$d_t = a$}
\put(6,3.25){$d_t = a,c$}
\put(6,0.25){$d_t = b$}
\put(12.25,10.5){\begin{minipage}{70pt}
$h_{t+1}=h_t-1$\\
$d_{t+1}=a,c$\\
$\phi_{t+1} \in \fpl_b$
\end{minipage}}
\put(12.25,7.5){\begin{minipage}{70pt}
$h_{t+1}=h_t$\\
$d_{t+1}=a$\\
$\phi_{t+1} \in \fpl_w$
\end{minipage}}
\put(12.25,4.5){\begin{minipage}{70pt}
$h_{t+1}=h_t+1$\\
$d_{t+1}=b,c$\\
$\phi_{t+1} \in \fpl_w$
\end{minipage}}
\put(12.25,1.5){\begin{minipage}{70pt}
$h_{t+1}=h_t$\\
$d_{t+1}=b$\\
$\phi_{t+1} \in \fpl_b$
\end{minipage}}
\end{picture}
\]
\caption{\label{fig.tabproof}Illustration of the analysis of Wieland
  half-gyration, involved in the proof of
  Lemma~\ref{lem.orbitpattern}. Graphical notations are as in the
  general treatment of Wieland gyration made in Section 3 of
  \cite{usRS} (see in particular \cite[Fig.~4]{usRS}, and Figure
  \ref{fig.LLpLm} here). More precisely, we denote in yellow only the
  plaquettes concerned with the operation $H$ (and not the ones
  concerned with $H^{-1}$), and we put a red bullet on those
  plaquettes that are left stable by gyration. An orange circle
  denotes the refinement position.  }
\end{figure}

As the orbit involves half-gyrations, given that $\phi_0 \in \fpl_+$,
$\phi_t \in \fpl_+$ if $t$ is even and $\phi_t \in \fpl_-$ if it is
odd. This implies that $\phi_t \in \fpl_b$ if $h_t-t$ is odd, and
$\phi_t \in \fpl_w$ if it is even (at time $t=0$, the black
terminations have labels $1,3,5,\ldots$). From the definition of $H$
in terms of $H_+$ and $H_-$, for both parities the black external
terminations are paired to the white termination \emph{at their
  right}, so that the black pattern rotates counter-clockwise, and the
white one clockwise. This implies that we have to study the local
behaviour of half-gyration only in the two situations described on the
left-most column of Figure~\ref{fig.tabproof} (the other two choices,
obtained by vertical reflection, would correspond to $H^{-1}$).

The two situations are very similar, and we discuss in detail only the
first case. Furthermore, in our graphical representation, we draw only
square plaquettes in a neighbourhood of the refinement position,
i.e.\ we describe the `generic' situation (in which we are far from
the corners, and no faces with less than 4 sides are present in the
neighbourhood). This is done only for simplicity of the visualisation,
and it is easily seen that the actual shape of the faces (within the
ones allowed for Wieland gyration), or the vicinity of corners, never
interfere with the local properties to be determined.

The well-known fact that there is a unique refinement position on a
reference side implies that, at any time and for both $H$ and
$H^{-1}$, there is at most one plaquette of the appropriate parity of
the form $\sqqV$ or $\sqqH$, and it is adjacent to the refinement
position (at its up-left or up-right corner, in our drawings). This
justifies the fact that we analyse only a neighbourhood of~$h_t$.

In studying $H$, for $\phi_t \in \fpl_b$, the two cases in which there
is one such plaquette, or there is not, occur exactly if $d_t \in
\{b,c\}$, or $d_t = a$, respectively, corresponding to the two
sub-cases of the first row in Figure \ref{fig.tabproof} (second
column).
For what concerns the neighbourhood of the reference side, the result
of the half-gyration is completely determined (and illustrated in the
third column of the drawing), from which we can deduce the properties
of $\phi_{t+1}$ summarised in the last column.  The collection of
these properties coincides with the statement of the lemma.  \qed


\bigskip
\noindent
We have thus determined that, in any orbit $\cO$,
$h_{t+1} - h_{t-1} \in \{-2,-1,+1,+2\}$. It is natural to say that the
configuration $\phi_t$ is \emph{ascending}, or \emph{descending}, if $h_{t+1} -
h_{t-1} \in \{+1,+2\}$, or $\{-1,-2\}$, respectively.  The table in
Lemma \ref{lem.orbitpattern} implies that a configuration $\phi$ is
descending if the edge incident to the refinement position is black,
and ascending if it is white.
As a consequence,
the sets $\cO_b^p$ and $\cO_w^p$ correspond to the sets of
configurations within the orbit that are descending or ascending,
respectively.

\begin{figure}[!tb]
\[
\includegraphics{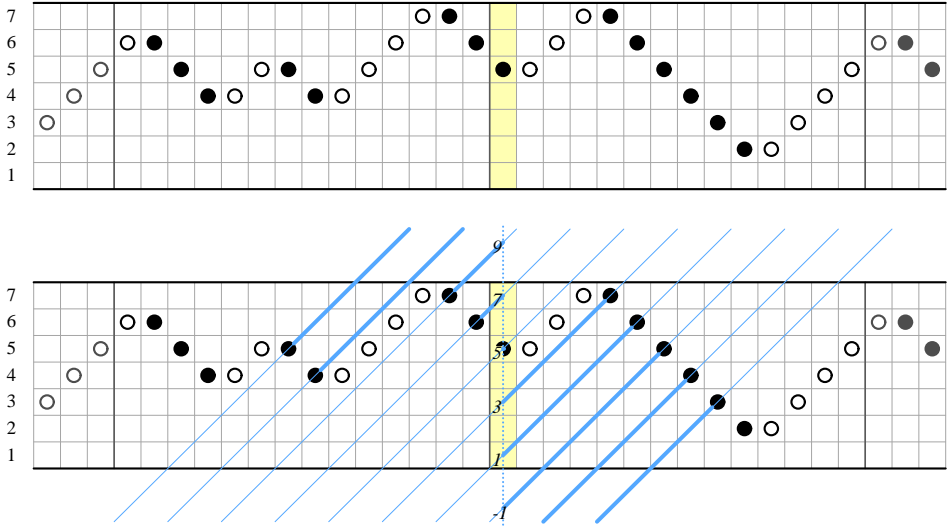}
\]
\caption{\label{fig.orbtraj}Top: a typical trajectory $h_t$, for
  $L=7$, and $p=28$. Time $t=0$ is at the yellow band.
  Bottom: illustration of Lemma \ref{h-t}, stating that, on the
  infinite orbit, for all odd values $c$ there exists a unique time
  $t$ such that $h_t-t=c$. Cyan diagonals are the level lines of $g_t
  = h_t-t$, at odd values.}
\end{figure}

The Lemma \ref{lem.orbitpattern} also implies
\begin{corollary}
\label{coroll.monotsubseq}
The sequences $h_t$ are composed of alternating ascending/desceding
monotonic subsequences of slope $\pm 1$, each of length at least 2.
The local minima have $d_t=a$, the local maxima $d_t=b$, and the other
elements have $d_t=c$.
In particular, local maxima and minima are achieved on plateaux of
length exactly 2.
\end{corollary}
\noindent
Some aspects of this corollary are illustrated in Figure
\ref{fig.orbtraj}, top, through an example.

The parity statement of Corollary \ref{coroll.monotsubseq} has an
important consequence.

\begin{lemma}
\label{coroll.monotsubs-equal-height}
Within a periodic portion $\cO^p$ of a gyration orbit $\cO$, the
number of FPL's that belong to $\cO^p_+$, $\cO^p_-$, $\cO^p_b$ and
$\cO^p_w$ at any given refinement position $h$ are all equal. In
in other words for any function $F$ on $\{1,\ldots,L\}$
and any orbit $\cO$ it holds
\be
\sum_{\phi\in \cO^p_+} F(h(\phi))
= \sum_{\phi\in \cO^p_-} F(h(\phi))
= \sum_{\phi\in \cO^p_{b}} F(h(\phi))
= \sum_{\phi\in \cO^p_{w}} F(h(\phi))
\ef.
\ee
\end{lemma}
\proof Fix an orbit $\cO$ and a value $h$.  Let $N_{\ast}(h)= |\{ \phi
\in \cO^p_{\ast},\; h(\phi)=h \}|$, for $\ast=+,-,b,w$. We want to
prove that $N_+(h)=N_-(h)=N_b(h)=N_w(h)$.  A consequence of Lemma
\ref{lem.orbitpattern} (the fact that $h_t-t$ has given parity for
configurations in $\fpl_{b/w}$) is that for $h$ odd $N_+(h)=N_b(h)$
and $N_-(h)=N_w(h)$, and for $h$ even $N_+(h)=N_w(h)$ and
$N_-(h)=N_b(h)$.
Call $\{ t_i \}_{i \in \mathbb{Z}}$ the ordered sequence of times
along the orbit such that $h(\phi_{t_i})=h$ (say, with $t_0$ the first
positive value).
Clearly $t_{i+N_b(h)+N_w(h)}=t_i + p$.
From the discrete continuity properties of $h_t$, the
configurations $\{ \phi_{t_i} \}_{i \in \mathbb{Z}}$ are alternating
ascending and descending, i.e., as a further consequence of
Lemma~\ref{lem.orbitpattern}, if $\phi_{t_i} \in \fpl_b$ then
$\phi_{t_{i\pm 1}} \in \fpl_w$, and vice versa.
This proves that $N_b(h)=N_w(h)$, and allows to
conclude. \qed


\bigskip
\noindent
Note that, as an outcome of the construction, we have natural
involutions between $\fpl_+$ and $\fpl_-$, and between $\fpl_b$ and
$\fpl_w$, which relate FPL's in the same orbit, preserve the
refinement position, and preserve the link pattern up to rotation, e.g.~by
associating to $\phi(t_{i}) \in \fpl_b$ the configuration 
$\phi(t_{i+1}) \in \fpl_w$.


The partition of $\fpl(\Lambda)$ into orbits is particularly useful
when one considers the image of a vector $\kket{v}\in
\mathbb{C}^{\fpl_+(\Lambda)}$ under $\Sym \circ \Pi_+$, or of a vector
$\kket{w}\in \mathbb{C}^{\fpl_b(\Lambda)}$ under $\Sym\circ\Pi_b$.
Indeed, writing $\kket{v}$ and $\kket{w}$ as
\begin{align}
\label{eq.563254876}
\kket{v} 
&= \sum_{\cO}  \sum_{\phi \in \cO^p_+} v(\phi)\kket{\phi}
\ef;
&
\kket{w} 
&= \sum_{\cO}  \sum_{\phi \in \cO^p_b} w(\phi)\kket{\phi}
\ef;
\end{align}
their images under
$\Sym \circ \Pi_+$ and $\Sym \circ
\Pi_b$ are simply 
\begin{align}
\label{eq.orbitdecompPM}
\Sym \ \Pi_+ \kket{v} 
&=
\sum_{\cO} 
\Big(\sum_{\phi \in \cO^p_+}
v(\phi)\Big) \ket{\cO}
\ef;
&
\Sym \ \Pi_b \kket{w} 
&=
\sum_{\cO} 
\Big(\sum_{\phi \in \cO^p_b}
w(\phi)\Big) \ket{\cO}
\ef;
\end{align}
where $\ket{\cO} =\Sym \ \Pi_+ \kket{\phi}$ for any $\phi\in \cO^p_+$,
or also $\Sym \ \Pi_b \kket{\phi}$ for any $\phi\in \cO^p_b$ (these
states are well-defined, i.e.~do not depend on the choice of
$\phi$, because
the link pattern is preserved by Wieland gyration, up to rotations).

\subsection{Proof of Di Francesco's 2004 conjecture}
\label{ssec.ph04proof}

\noindent
Recall that we defined
\be
\ket{\Psi_\Lambda(t)} = \sum_{\phi \in \fpl_b(\Lambda)}
t^{h(\phi)-1} \ket{\Pi_b(\phi)}
\ef;
\ee
as the refinement of the enumerations according to the `new'
function $\tilde{\pi}(\phi) = \Pi_b(\phi)$.
Then let 
\be
\ket{\Psi'_\Lambda(t)} = \sum_{\phi \in \fpl_+(\Lambda)}
t^{h(\phi)-1} \ket{\Pi_+(\phi)}
\ef;
\ee
the refinement of the enumerations according to the `ordinary'
function $\pi(\phi) = \Pi_+(\phi)$.

At this point, and at the light of Theorem \ref{main}, it is easy to
restate and prove Di~Francesco's conjecture of \cite{PdF04}.
\begin{theorem}[Di Francesco's 2004 conjecture]
For any dihedral domain $\Lambda$ and reference side $r$,
\label{df-conj}
\be
\label{df-eqs}
\Sym\, \ket{\Psi'_\Lambda(t)} = \Sym\, \ket{\Psi_\Lambda(t)}
\ef.
\ee 
\end{theorem}
\vspace{1mm}
\noindent
More precisely, the original conjecture states that, for $\Lambda$ the
$n \times n$ square domain,
\be
\Sym\, \ket{\Psi'_{\Lambda}(t)}
=
\Sym\, 
\ket{\Psi^{(1)}_{O(1)}(t)}
\ef.
\ee
However, this fact naturally extends to all dihedral domains, up to a
proportionality factor, namely
$
\Sym\, \ket{\Psi'_{\Lambda}(t)}
=
K_{\Lambda}(t)\;
\Sym\, 
\ket{\Psi^{(1)}_{O(1)}(t)}
$, and, at the light of Theorem \ref{main} (in the formulation of 
equation (\ref{eq.KinTh})), the restatement in Theorem \ref{df-conj} follows.

\bigskip
\proofof{Theorem \ref{df-conj}}
Consider the vectors
$\kket{s'_\Lambda(t)} \in \mathbb{C}^{\fpl_+(\Lambda)}$ and
$\kket{s_\Lambda(t)} \in \mathbb{C}^{\fpl_b(\Lambda)}$,
\begin{align}
\kket{s'_\Lambda(t)}
&=\sum_{\phi \in \fpl_+(\Lambda)}
t^{h(\phi)-1}\kket{\phi}
\ef;
&
\kket{s_\Lambda(t)}
&=\sum_{\phi \in \fpl_b(\Lambda)}
t^{h(\phi)-1}\kket{\phi}
\ef.
\end{align}
Thus $\ket{\Psi'_\Lambda(t)}= \Pi_+ \kket{s'_\Lambda(t)}$ and
$\ket{\Psi_\Lambda(t)}= \Pi_b \kket{s_\Lambda(t)}$.  In view of
equations (\ref{eq.orbitdecompPM}), and of the decompositions
\begin{align}
\kket{s'_\Lambda(t)}
&=
\sum_{\cO} \sum_{\phi \in \cO^p_+} 
t^{h(\phi)-1} \kket{\phi}
\ef;
&
\kket{s_\Lambda(t)}
&=
\sum_{\cO} \sum_{\phi \in \cO^p_b} 
t^{h(\phi)-1} \kket{\phi}
\ef;
\end{align}
(so that we are in the situation of equations (\ref{eq.563254876})),
a sufficient condition for the theorem to hold is that, for all the
orbits $\cO$,
\be
\sum_{\phi \in \cO^p_+}
t^{h(\phi)-1}=\sum_{\phi \in \cO^p_b}
t^{h(\phi)-1}
\ef,
\ee
which is a special case of Lemma
\ref{coroll.monotsubs-equal-height}, with $F(x) = t^{x-1}$.
\qed


\bigskip
\noindent
As a corollary of equation (\ref{df-eqs}) we have the ordinary
Razumov--Stroganov correspondence
\begin{corollary}[Ordinary Razumov--Stroganov correspondence]
\label{usualRS1} 
The state $\ket{\Psi'_\Lambda(1)}$ satisfies equation
(\ref{RS-eq0}).
\end{corollary}
\proof The Wieland Theorem shows that the state
$\ket{\Psi'_\Lambda(1)}$ is rotationally invariant.  From the
reasonings in Section \ref{exch-section} for
$\ket{\Psi^{(1)}_{O(1)}(1)}$, and the proportionality of
$\ket{\Psi^{(1)}_{O(1)}(1)}$ and $\ket{\Psi_\Lambda(1)}$ stated by
Theorem \ref{main}, we already know that $\ket{\Psi_\Lambda(1)}$ is
the unique solution of equation (\ref{RS-eq0}) up to normalisation,
and in particular it is rotationally invariant.
Thus equation (\ref{df-eqs}) is equivalent to
$\ket{\Psi'_\Lambda(1)}=\ket{\Psi_\Lambda(1)}$, with no need of
symmetrisation.  As a result, $\ket{\Psi'_\Lambda(1)}$ satisfies
equation (\ref{RS-eq0}).  \qed

\subsection{A bijection between $\fpl_+(\Lambda)$ and
  $\fpl_b(\Lambda)$}
\label{bij-section} 

\noindent
In the previous paragraphs we have introduced bijections between
$\fpl_+$ and $\fpl_-$, and between $\fpl_b$ and $\fpl_w$, which relate
FPL's in the same orbit, preserve the refinement position, and
preserve the link pattern up to rotation.

In this section we introduce an explicit bijection between
$\fpl_+(\Lambda)$ and $\fpl_b(\Lambda)$, which relates FPL's within
the same orbit, and preserves the link pattern, when produced with
$\Pi_+$ and $\Pi_b$ respectively (but does \emph{not} preserve the
refinement position).

This bijection allows to recover the ordinary Razumov--Stroganov
correspondence based on fully-packed loop configurations in
$\fpl_+(\Lambda)$, from the central theorem proven in this paper,
Theorem \ref{main}, without using Theorem \ref{df-conj}, nor the
results in \cite{pdf-pzj1}.

Roughly speaking, the idea is to rotate a FPL $\phi\in\fpl_+(\Lambda)$
until the image of the external edge with label $1$ coincides with the
refinement position. While \emph{a priori} it is not obvious that this
event ever occurs within the orbit, this fact is ensured by the
following
\begin{lemma}
\label{h-t}
The function $g_t=h_t-t$ is non-increasing, and each odd value has exactly one
preimage.
\end{lemma}
\proof The proof is a simple consequence of the table in Lemma
\ref{lem.orbitpattern}, which implies that if $g_t$ is even, then
$g_{t+1} \in \{g_t, g_t -1\}$, and if $g_t$ is odd, then $g_{t+1} \in
\{g_t -2, g_t -1\}$.  The idea is also illustrated in Figure
\ref{fig.orbtraj}, bottom.
\qed

\bigskip
\noindent
This lemma provides the claimed bijection.  For the orbit with
$\phi_0=\phi$, call $t^*(\phi)$ the preimage of $1$ under $g_t$ (which
is unique from the lemma). Then we have
\begin{proposition}
We have a bijection
$\Theta:\fpl_+(\Lambda)\rightarrow\fpl_b(\Lambda)$, defined, together
with its inverse, as
\begin{align}
\label{eq.654765}
\Theta(\phi) 
&:= H^{t^*(\phi)} \phi
\ef;
&
\Theta^{-1}(\phi) 
&:= H^{-h(\phi)+1} \phi
\ef.
\end{align}
Furthermore, $\Pi_+(\phi) =\Pi_b(\Theta(\phi))$. 
\end{proposition}
\proof
The fact that $\Pi_+(\phi) =\Pi_b(\Theta(\phi))$ is a restatement of
$g_t=1$ (at the light of the fact that each half-gyration rotates the
link pattern one position counter-clockwise).


The fact that $\Theta(\Theta^{-1}(\phi))=\phi$, as defined in
(\ref{eq.654765}), follows from the definition of $t^*(\phi)$.  The
fact that also $\Theta^{-1}(\Theta(\phi))=\phi$,
follows by observing that $0= h_{t^*(\phi)} - t^*(\phi) -1 =
h(\Theta(\phi)) - t^*(\phi) -1$, that is, $(H^{t^*(\phi)})^{-1} =
H^{-(h(\Theta(\phi))-1)}$.
\qed




\bigskip
\noindent
Using the statement
$\Pi_+(\phi) =\Pi_b(\Theta(\phi))$, we can easily show that
$\ket{\Psi'_\Lambda(1)}= 
\ket{\Psi_\Lambda(1)}$:
\be
\ket{\Psi'_\Lambda(1)} 
= \Pi_+ \!\!\sum_{\phi' \in \fpl_+(\Lambda)}\!\! \kket{\phi'} 
= \Pi_b \!\!\sum_{\phi' \in \fpl_+(\Lambda)}\!\! \kket{\Theta(\phi')} 
= \Pi_b \!\!\sum_{\phi \in \fpl_b(\Lambda)}\!\! \kket{\phi}
= \ket{\Psi_\Lambda(1)}
\ef.
\ee
Through Theorem \ref{main}, this is another proof of the
Razumov--Stroganov correspondence.


\end{document}